\documentclass{amsart}
\usepackage{amssymb}
\usepackage{amsmath}
\usepackage{amsfonts}

\newtheorem{theorem}{Theorem}
\theoremstyle{plain}

\newtheorem{lemma}{Lemma}

\newtheorem{proposition}{Proposition}
\newtheorem{remark}{Remark}

\numberwithin{equation}{section}

\begin{document}
\title[On the first eigenvalue for a $(p(x),q(x))$-Laplacian elliptic system]%
{On the first eigenvalue for a $(p(x),q(x))$-Laplacian elliptic system }
\subjclass[2010]{35J60, 35P30, 47J10, 35A15, 35D30}
\keywords{$p(x)$-Laplacian, Variable exponents, Weak solution, Eigenvalue, Regularity, Boundedness.}

\begin{abstract}
In this article, we deal about the first eigenvalue for a nonlinear gradient type elliptic system involving variable exponents growth conditions. Positivity, boundedness and regularity of associated eigenfunctions for auxiliaries systems are established.
\end{abstract}

\author{Abdelkrim Moussaoui}
\address{Abdelkrim Moussaoui\\
Applied Mathematics Laboratory (LMA), Faculty of Exact Sciences,\\
A. Mira Bejaia University, Targa Ouzemour 06000 Bejaia, Algeria}
\email{abdelkrim.moussaoui@univ-bejaia.dz}
\author{Jean V\'{e}lin}
\address{Jean V\'{e}lin\\
D\'{e}partement de Math\'{e}matiques et Informatique, Laboratoire CEREGMIA,
Universit\'{e} des Antilles et de La Guyane, Campus de Fouillole 97159
Pointe-\`{a}-Pitre, Guadeloupe (FWI)}
\email{jean.velin@univ-ag.fr}
\maketitle

\section{Introduction and setting of the problem}

\label{S1}

In the present paper, we focus on finding a non zero first eigenvalue for
the system of quasilinear elliptic equations%
\begin{equation}
(P)\,\,\left\{ 
\begin{array}{ll}
-\Delta _{p(x)}{u}=\lambda c(x)(\alpha (x)+1)|{u}|^{\alpha (x)-1}u|{v}%
|^{\beta (x)+1} & \hbox{ in }\,\,\Omega \\ 
-\Delta _{q(x)}{v}=\lambda c(x)(\beta (x)+1)|{u}|^{\alpha (x)+1}{v}|{v}%
|^{\beta (x)-1} & \hbox{ in }\,\,\Omega \\ 
{u}={v}=0 & \hbox{ on }\,\,\partial \Omega%
\end{array}%
\right.  \label{p}
\end{equation}%
on a bounded domain $\Omega \subset \mathbb{R}^{N}.$ Here $\Delta
_{p(x)}u=div(|\nabla u|^{p(x)-2}\nabla u)$ and $\Delta _{q(x)}v=div(|\nabla
v|^{q(x)-2}\nabla u)$ are usually named the $p(x)$-Laplacian and the $q(x)$%
-Laplacian operator. \newline
During the last decade, the interest for partial differential equations
involving the $p(x)$-Laplacian operator is increasing. When the exponent
variable function $p(\cdot )$ is reduced to be a constant, $\Delta _{p(x)}u$
becomes the well-known $p$-Laplacian operator $\Delta _{p}u.$ The $p(x)$%
-Laplacian operator possesses more complicated nonlinearity than the $p$%
-Laplacian. So, one cannot always to transpose to the problems arising the $%
p(x)$-Laplacian operator the results obtained with the $p$-Laplacian. The
treatments of solving these problem are often very complicated and needs a
mathematical tools (Lebesgue and Sobolev spaces with variable exponents, see
for instance \cite{dhhr} and its abundant reference). Among them, finding
first eigenvalue of $p(x)$-Laplacian Dirichlet presents more singular
phenomena which do not appear in the constant case. More precisely, it is
well known that the first eigenvalue for the $p(x)$-Laplacian Dirichlet
problem may be equal to zero (for details, the reader interested can consult 
\cite{fanandall}). In \cite{fanandall}, the authors consider that $\Omega $
is a bounded domain and $p$ is a continuous function from ${\overline{\Omega 
}}$ to $]1,+\infty \lbrack .$ They given some geometrical conditions
insuring that the first eigenvalue is $0.$ Otherwise, in one dimensional
space, monotonicity assumptions on the function $p$ is a necessary and
sufficient condition such that the first eigenvalue is strictly positive. In
higher dimensional case, assuming monotonicity of an associated function
defined by $p,$ the first eigenvalue is strictly positive.

The fact of the first eigenvalue is zero, has been observed earliest by \cite%
{fz}. Indeed, the authors illustrate this phenomena by taking $\Omega
=(-2,2) $ and $p(x)=3\chi _{\lbrack 0,1]}(x)+(4-|x|)\chi _{\lbrack 1,2]}(x).$
In this condition, the Rayleigh quotient 
\begin{equation*}
\begin{array}{l}
\mu _{1}=\inf_{u\in W_{0}^{1,p(\cdot )}(\Omega )\setminus \{0\}}\frac{%
\int_{\Omega }|\nabla u|^{p(x)}}{\int_{\Omega }|u|^{p(x)}}%
\end{array}%
\end{equation*}
is equal to zero. The main reason derives that the well-known Poincar\'{e}
inequality is not always fulfilled. However, Fu in \cite{fu} shown that when 
$\Omega $ is a bounded Lipschitz domain, $p$ is $L^{\infty }(\Omega )$ the
Poincar\'{e}{\ }inequality holds (i.e. there is a constant $C$ depending on $%
\Omega $ such that for any $u\in W_{0}^{1,p(x)}(\Omega ),$ $\int_{\Omega
}|u|^{p(x)}\leq C\int_{\Omega }|\nabla u|^{p(x)}$). For a use of this result
see for instance \cite{cf}, \cite{jf}.

Further works established suitable conditions drawing to a non zero first
eigenvalue (see \cite{f}, \cite{mr}, \cite{fm}).

Compared the investigation for one equation, elliptic systems haven't a
similar growth concerning in the first eigenvalue. First of all, when $p(x)$
and $q(x)$ are constant on $\Omega ,$ in \cite{dt}, the following elliptic
Dirichlet system is considered 
\begin{equation}  \label{Spq}
\left\{ 
\begin{array}{ll}
-\Delta _{p}u=\lambda u|u|^{\alpha -1}|v|^{\beta +1} & \hbox{ in }\,\,\Omega
\\ 
-\Delta _{q}v=\lambda |u|^{\alpha +1}|v|^{\beta +1}v & \hbox{ in }\,\,\Omega
\\ 
u=v=0 & \hbox{ on }\,\,\partial \Omega .%
\end{array}%
\right.
\end{equation}%
Assuming, $\Omega $ is a bounded open in $%
\mathbb{R}
^{N}$ with smooth boundary $\partial \Omega $ and the constant exponents $%
-1<\alpha ,\beta $ and $1<p,q<N$ satisfying the condition 
\begin{equation}
\begin{array}{l}
C_{p,q}^{\alpha ,\beta }:\,\frac{\alpha +1}{p}+\frac{\beta +1}{q}=1\,\,%
\hbox{ and }(\alpha +1)\frac{N-p}{Np}+(\beta +1)\frac{N-q}{Nq}<1,%
\end{array}
\label{Cpq}
\end{equation}%
the author shown the existence of the first eigenvalue $\lambda (p,q)>0$
associated to a positive and unique eigenfunction $(u^{\ast },v^{\ast }).$
Further more, this result have been extended by Kandilakis and al. \cite{kmz}
for the system 
\begin{equation}
\left\{ 
\begin{array}{ll}
\Delta _{p}u+\lambda a(x)|u|^{p-2}u+\lambda b(x)u|u|^{\alpha -1}|v|^{\beta
+1}=0 & \hbox{ in }\,\,\Omega \\ 
\Delta _{q}v+\lambda d(x)|u|^{p-2}u+\lambda b(x)|u|^{\alpha +1}|v|^{\beta
+1}v=0 & \hbox{ on }\,\,\Omega \\ 
|\nabla u|^{p-2}\nabla u\cdot \nu +c_{1}(x)|u|^{p-2}u=0 & \hbox{ on }%
\,\,\partial \Omega \\ 
|\nabla v|^{q-2}\nabla v\cdot \nu +c_{2}(x)|v|^{q-2}v=0 & \hbox{ on }%
\,\,\partial \Omega ,%
\end{array}%
\right.  \label{Spqabcd}
\end{equation}%
with $\Omega $ is an unbounded domain in $%
\mathbb{R}
^{N}$ with non compact and smooth boundary $\partial \Omega ,$ the constant
exponents $0<\alpha ,\beta $ and $1<p,q<N$ satisfying 
\begin{equation}
\begin{array}{l}
\frac{\alpha +1}{p}+\frac{\beta +1}{q}=1\,\,\hbox{
and }(\alpha +1)\frac{N-p}{Np}<q,\,\,(\beta +1)\frac{N-q}{Nq}<p.%
\end{array}
\label{Cpqabcd}
\end{equation}
Inspired by \cite{dt}, Khalil and al. in \cite{kmo} shown that the first
eigenvalue $\lambda _{p,q}$ of (\ref{Spq}) is simple and moreover they
established stability (continuity) for the function $(p,q)\longmapsto
\lambda (p,q).$

Motivated by the aforementioned papers, in this work we establish the
existence of one-parameter family of nontrivial solutions $((\hat{u}_{R},%
\hat{v}_{R}),\lambda _{R}^{\ast })$ for all $R>0$ for problem (\ref{p}). In
addition, we show that the corresponding eigenfunction $(\hat{u}_{R},\hat{v}%
_{R})$ is positive in $\Omega $, bounded in $L^{\infty }(\Omega )\times
L^{\infty }(\Omega )$ and belongs to $C^{1,\gamma }(\overline{\Omega }%
)\times C^{1,\gamma }(\overline{\Omega })$ for certain $\gamma \in (0,1)$ if 
$p,q\in C^{1}(\overline{\Omega })\cap C^{0,\theta }(\overline{\Omega })$.
Furthermore, by means of geometrical conditions on the domain $\Omega $, we
prove that the infimum of the eigenvalues of (\ref{p}) is positive. To the
best of our knowledge, it is for the first time when the positive infimum
eigenvalue for systems involving $p(x)$-Laplacian operator is studied.
However, we point out that in this paper, the existence of an eigenfunction
corresponding to the infimum of the eigenvalues of (\ref{p}) is not
established and therefore, this issue still remains an open problem.

The rest of the paper is organized as follows. Section \ref{S2} contains
hypotheses, some auxiliary and useful results involving variable exponent
Lebesgue-Sobolev spaces and our main results. Section \ref{S3} and section %
\ref{S4} present the proof of our main results.

\section{ Hypotheses - Main results and some auxiliary results}

\label{S2}

Let $L^{p(x)}(\Omega )$ be the generalized Lebesgue space that consists of
all measurable real-valued functions $u$ satisfying%
\begin{equation*}
\begin{array}{l}
\rho _{p(x)}(u)=\int_{\Omega }|u(x)|^{p(x)}dx<+\infty ,%
\end{array}%
\end{equation*}
endowed with the Luxemburg norm%
\begin{equation*}
\begin{array}{l}
\left\Vert u\right\Vert _{p(x)}=\inf \{\tau >0:\rho _{p(x)}(\frac{u}{\tau }%
)\leq 1\}.%
\end{array}%
\end{equation*}%
The variable exponent Sobolev space $W_{0}^{1,p(\cdot )}(\Omega )$ is
defined by%
\begin{equation*}
\begin{array}{l}
W_{0}^{1,p(x)}(\Omega )=\{u\in L^{p(x)}(\Omega ):|\nabla u|\in
L^{p(x)}(\Omega )\}.%
\end{array}%
\end{equation*}
The norm $\left\Vert u\right\Vert _{1,p(x)}=\left\Vert \nabla u\right\Vert
_{p(x)}$ makes $W_{0}^{1,p(x)}(\Omega )$ a Banach space and the following
embedding%
\begin{equation}
\begin{array}{l}
W_{0}^{1,p(x)}\hookrightarrow L^{r(x)}(\Omega )%
\end{array}
\label{20}
\end{equation}%
is compact with $1<r(x)<\frac{Np(x)}{N-p(x)}$. 

\subsection{ Hypotheses}

\begin{description}
\item[\textrm{(H.1)}] $\Omega $ $\hbox{ is an bounded open of  }\mathbb{R}%
^{N},$ $\hbox{ its boundary }\partial \Omega \hbox{ of class }C^{2,\delta },$
for certain $0<\delta <1,$

\item[\textrm{(H.2)}] $c:\Omega \longrightarrow \mathbb{R}_{+}$ and$\ c\in
L^{\infty }(\Omega ),$

\item[\textrm{(H.3)}] $\alpha ,\beta :\overline{\Omega }\rightarrow
]1,+\infty \lbrack \hbox{ two continuous functions satisfying  }$

\begin{equation*}
1<\alpha ^{-}=\inf_{x\in \Omega }\alpha (x)\leq \alpha ^{+}=\sup_{x\in
\Omega }\alpha (x)<\infty , \newline
1<\beta ^{-}=\inf_{x\in \Omega }\beta (x)\leq \beta ^{+}=\sup_{x\in \Omega
}\beta (x)<\infty
\end{equation*}
and%
\begin{equation*}
\begin{array}{l}
\frac{\alpha (x)+1}{p(x)}+\frac{\beta (x)+1}{q(x)}=1,%
\end{array}%
\end{equation*}

\item[\textrm{(H.4)}] $p\hbox{ and }q%
\hbox{ are two variable exponents of
class }C^{1}(\overline{\Omega })\hbox{ satisfying }$

\begin{equation*}
p(x)<\frac{Np(x)}{N-p(x)}, q(x)<\frac{Nq(x)}{N-q(x)},
\end{equation*}%
with%
\begin{equation*}
\begin{array}{l}
1<p^{-}=\inf_{x\in \Omega }p(x)\leq p^{+}=\sup_{x\in \Omega }p(x)<\infty ,
\\ 
1<q^{-}=\inf_{x\in \Omega }q(x)\leq q^{+}=\sup_{x\in \Omega }q(x)<\infty .%
\end{array}%
\end{equation*}
\end{description}


\subsection{ Main results}

\medskip Throughout this paper, we set $X_{0}^{p(x),q(x)}(\Omega
)=W_{0}^{1,p(x)}(\Omega )\times W_{0}^{1,q(x)}(\Omega ).$

Define on $X_{0}^{p(x),q(x)}(\Omega )$ the functionals $\mathcal{A}$ and $%
\mathcal{B}$ as follows: 
\begin{equation}
\begin{array}{l}
\mathcal{A}(z,w)=\int_{\Omega }\frac{1}{p(x)}\left\vert \nabla z\right\vert
^{p(x)}dx+\int_{\Omega }\frac{1}{q(x)}\left\vert \nabla w\right\vert
^{q(x)}dx,%
\end{array}
\label{A}
\end{equation}%
\begin{equation}
\begin{array}{l}
\mathcal{B}(z,w)=\int_{\Omega }c(x)|z|^{\alpha (x)+1}|w|^{\beta (x)+1}dx,%
\end{array}
\label{B}
\end{equation}%
and denote by $\left\Vert (z,w)\right\Vert =\left\Vert z\right\Vert
_{1,p(x)}+\left\Vert w\right\Vert _{1,q(x)}$. The same reasoning exploited
in \cite{FZ} implies that $\mathcal{A}$ and $\mathcal{B}$ are of class $%
C^{1}(X_{0}^{p(x),q(x)}(\Omega ),%
\mathbb{R}
).$ The Fr\'{e}chet derivatives of $\mathcal{A}$ and $\mathcal{B}$ at $(z,w)$
in $X_{0}^{p(x),q(x)}(\Omega )$ are given by%
\begin{equation}
\begin{array}{l}
\mathcal{A}^{\prime }(z,w)\cdot (\varphi ,\psi )=\int_{\Omega }|\nabla {z}%
|^{p(x)-2}\nabla {z}\cdot \nabla {\varphi } {\ }dx+\int_{\Omega }|\nabla {w}%
|^{q(x)-2}\nabla {w\cdot \nabla \psi }{\ }dx%
\end{array}
\label{Aprime}
\end{equation}

and%
\begin{equation}
\begin{array}{c}
\mathcal{B}^{\prime }(z,w)\cdot (\varphi ,\psi )=\int_{\Omega }c(x)(\alpha
(x)+1)|z|^{\alpha (x)-1}|w|^{\beta (x)+1}\varphi \\ 
+\int_{\Omega }c(x)(\beta (x)+1)|z|^{\alpha (x)+1}|w|^{\beta (x)-1}{w\psi }{%
\ }dx,%
\end{array}
\label{Bprime}
\end{equation}%
where $(\varphi ,\psi )\in X_{0}^{p(x),q(x)}(\Omega ).$

Let $R>0$ be fixed, we set%
\begin{equation*}
\mathcal{X}_{R}=\{(z,w)\in X_{0}^{p(x),q(x)}(\Omega );\,\,\mathcal{B}%
(z,w)=R\}.
\end{equation*}
It is obvious to notice that the set $\mathcal{X}_{R}$ is not empty. Indeed,
let $(z_{0},w_{0})\in X_{0}^{p(x),q(x)}(\Omega )$ such that $\mathcal{B}%
(z_{0},w_{0})=b_{0}>0,$ if $b_{0}=R,$ we done. Otherwise, for $z_{R}=({R}/{%
b_{0}})^{1/p(x)}z_{0}$ and $w_{R}=({R}/{b_{0}})^{1/q(x)}w_{0},$ it is easy
to note that $\mathcal{B}(z_{R},w_{R})=R.$

Now, define the Rayleigh quotients 
\begin{equation}
\begin{array}{l}
\lambda _{R}^{\ast }=\inf_{(z,w)\in \mathcal{X}_{R}}\frac{\mathcal{A}(z,w)}{%
\mathcal{B}(z,w)},%
\end{array}
\label{c1}
\end{equation}%
\begin{equation*}
\begin{array}{l}
\lambda _{p(x),q(x)}^{\ast }=\inf_{(z,w)\in X_{0}^{1,p(x),q(x)}(\Omega
)\setminus \{0\}}\frac{\mathcal{A}(z,w)}{\mathcal{B}(z,w)}%
\end{array}%
\end{equation*}
and%
\begin{equation}
\begin{array}{l}
{\lambda _{\ast }}_{R}=\inf_{(z,w)\in \mathcal{X}_{R}}\frac{\mathcal{A}(z,w)%
}{\int_{\Omega }c(x)(\alpha (x)+\beta (x)+2)|z|^{\alpha (x)+1}|w|^{\beta
(x)+1}dx}.%
\end{array}
\label{c2}
\end{equation}

\begin{remark}
\label{R1}The constant $\lambda _{R}^{\ast }$ in (\ref{c1}) can be written
as follows%
\begin{equation}
\begin{array}{l}
R\lambda _{R}^{\ast }=\inf_{\{(z,w)\in \mathcal{X}_{R}\}}\mathcal{A}(z,w).%
\end{array}
\label{infAR}
\end{equation}
\end{remark}

\medskip \medskip Our first main result provides the existence of a one -
parameter family of solutions for the system (\ref{p}).

\begin{theorem}
\label{T0} Assume that \textrm{(H.1)} - \textrm{(H.4)} hold. Then, the
system (\ref{p}) has a one-parameter family of nontrivial solutions $\left( (%
\hat{u}_{R},\hat{v}_{R}),\lambda _{R}^{\ast }\right) $ for all $R\in
(0,+\infty ).$ Moreover, if one of the following conditions holds:

\begin{description}
\item[\textrm{(a.1)}] There is vectors $l_{1},l_{2}\in \mathbb{R}%
^{N}\setminus \{0\}$ such that for all $x\in \Omega ,$ $%
f(t_{1})=p(x+t_{1}l_{1})$ and $g(t_{2})=q(x+t_{2}l_{2})$ are monotone for $%
t_{i}\in I_{i,x}=\{t_{i};\,\,x+t_{i}l_{i}\in \Omega \},$ $i=1,2.$

\item[\textrm{(a.2)}] There is $x_{1},x_{2}\notin \overline{\Omega }$ such
that for all $w_{1},w_{2}\in \mathbb{R}\setminus \{0\}$ with $\Vert
w_{1}\Vert ,\Vert w_{2}\Vert =1,$ the functions $%
f(t_{1})=p(x_{0}+t_{1}w_{1}) $ and $g(t_{2})=p(x_{2}+t_{2}w_{2})$ are
monotone for $t_{i}\in I_{x_{i},w_{i}}=\{t_{i}\in \mathbb{R}%
;\,\,x_{i}+t_{i}w_{i}\in \Omega \},$ $i=1,2.$
\end{description}

Then, $\lambda _{p(x),q(x)}^{\ast }=\inf_{R>0}\lambda _{R}^{\ast }>0$ is the
positive infimum eigenvalue of problem (\ref{p}).
\end{theorem}

A second main result consists in positivity, boundedness and regularity for
the obtained solution of problem (\ref{p}).

\begin{theorem}
\label{T2}Let $R$ be a fixed and strictly positive real. Assume that \textrm{%
(H.3)} holds.\newline
Then, $(\hat{u}_{R},\hat{v}_{R})$ the nontrivial solution of problem (\ref{p}%
) is positive and bounded in $L^{\infty }(\Omega )\times L^{\infty }(\Omega
).$ Moreover, if $p,q\in C^{1}(\overline{\Omega })\cap C^{0,\gamma }(%
\overline{\Omega })$ for certain $\gamma \in (0,1)$ then $(\hat{u}_{R},\hat{v%
}_{R})$ belongs to $C^{1,\delta }(\overline{\Omega })\times C^{1,\delta }(%
\overline{\Omega }),$ $\delta \in (0,1)$.
\end{theorem}

The proof of Theorem \ref{T0} will be done in section \ref{S3} while in
section \ref{S4} we will present the proof of Theorem \ref{T2}.


\subsection{Some preliminaries lemmas}

\begin{lemma}
\label{L1}$(i)$ For any $u\in L^{p(x)}(\Omega )$ we have%
\begin{equation*}
\begin{array}{l}
\left\Vert u\right\Vert _{p(x)}^{p^{-}}\leq \rho _{p(x)}(u)\leq \left\Vert
u\right\Vert _{p(x)}^{p^{+}}\hbox{ \ if \ }\left\Vert u\right\Vert _{p(x)}>1,%
\end{array}%
\end{equation*}

\begin{equation*}
\begin{array}{l}
\left\Vert u\right\Vert _{p(x)}^{p^{+}}\leq \rho _{p(x)}(u)\leq \left\Vert
u\right\Vert _{p(x)}^{p^{-}}\hbox{ \ if \ }\left\Vert u\right\Vert
_{p(x)}\leq 1.%
\end{array}%
\end{equation*}

$(ii)$ For $u\in L^{p(x)}(\Omega )\backslash \{0\}$ we have 
\begin{equation}
\left\Vert u\right\Vert _{p(x)}=a\hbox{ \ if and only if }\rho _{p(x)}(\frac{%
u}{a})=1.  \label{normro}
\end{equation}
\end{lemma}

\begin{lemma}[{\protect\cite[Theorem 8.2.4]{dhhr}}]
\label{L2}For every $u\in W_{0}^{1,p(\cdot )}(\Omega )$ the inequality%
\begin{equation}
\Vert u\Vert _{L^{p(\cdot )}(\Omega )}\leq C_{N,p}\Vert \nabla u\Vert
_{L^{p(\cdot )}(\Omega )},  \label{unabla}
\end{equation}%
holds with a constant $C_{N,p}>0$.
\end{lemma}

Recall that if there exist a constant $L>0$ and an exponent $\theta \in
(0,1) $ such that 
\begin{equation*}
|p(x_{1})-p(x_{2})|\leq L|x_{1}-x_{2}|^{\theta }\hbox{ \ for all }%
x_{1},x_{2}\in \overline{\Omega },
\end{equation*}
then the function $p$ is said to be H\"{o}lder continuous on $\overline{%
\Omega }$ and we denote $p\in C^{0,\theta }(\overline{\Omega })$.

For a later use, we have the next result.

\begin{lemma}
\label{L6}For $s\in (0,1)$ it holds%
\begin{equation*}
\begin{array}{l}
{\sum }_{n=1}^{r}(n-1)s^{n-1}\leq \frac{s}{(s-1)^{2}}.%
\end{array}%
\end{equation*}
\end{lemma}

\textbf{Proof.} \ignorespaces
Recall that for $s>0$ we have%
\begin{equation*}
\begin{array}{l}
s^{r}-1=(s-1)(s^{r-1}+s^{r-2}+...+s+1),\forall r\in 
\mathbb{N}
^{\ast }.%
\end{array}%
\end{equation*}%
Multiplying by $s$ one get%
\begin{equation*}
\begin{array}{l}
rs^{r}=(s^{r}+s^{r-1}+...+s)+(s-1)((r-1)s^{r-1}+...+s) \\ 
=\frac{s-s^{r+1}}{1-s}+(s-1)((r-1)s^{r-1}+...+s),%
\end{array}%
\end{equation*}
for all $s\neq 1$. Thus, it follows that%
\begin{equation*}
\begin{array}{l}
(r-1)s^{r-1}+...+s=\frac{s-s^{r+1}}{(s-1)^{2}}-\frac{rs^{r-1}}{1-s}.%
\end{array}%
\end{equation*}
Hence, for $0<s<1$ one has%
\begin{equation*}
\begin{array}{l}
(r-1)s^{r-1}+...+s\leq \frac{s}{(s-1)^{2}}.%
\end{array}%
\end{equation*}
{\ \rule{0.2cm}{0.2cm}{\ }}

\section{Proof of Theorem \protect\ref{T0}}

\label{S3}

Taking account of the assumption \textrm{(H.3)}, we note that the system (%
\ref{p}) is arising from a nonlinear eigenvalue type problem. Solvability of
general class of nonlinear eigenvalues problems of type $\mathcal{A}^{\prime
}(x)=\lambda \mathcal{B}^{\prime }(x)$ have been treated by M.S Berger in 
\cite{msb}. We recall this main tool.

\begin{theorem}
\cite{msb} Suppose that the $C^1$ functionals $\mathcal{A}$ and $\mathcal{B}$
defined on the reflexive Banach space $X$ have the following properties:

\begin{enumerate}
\item $\mathcal{A}$ is weakly lower semicontinuous and coercive on $X\cap \{%
\mathcal{B}(x)\leq const. \};$

\item $\mathcal{B}$ is continuous with respect to weak sequential
convergence and $\mathcal{B}^{\prime}(x)=0$ only at $x=0.$
\end{enumerate}

Then the equation $\mathcal{A}^{\prime}(x)=\lambda\mathcal{B}^{\prime}(x)$
has a one-parameter family of nontrivial solutions $(x_R,\lambda_R)$ for all 
$R$ in the range of $\mathcal{B}(x)$ such that $\mathcal{B}(x_R)=R;$ and $%
x_R $ is characterized as the minimum of $\mathcal{A}(x)$ over the set $\{%
\mathcal{B}(x)=R\}.$ \label{msb}
\end{theorem}

\begin{remark}
In the statement (\textbf{ii}) of the theorem \ref{msb}, the condition
\textquotedblleft $\mathcal{B}^{\prime }(x)=0$ only at $x=0$%
\textquotedblright \thinspace may be replaced by \textquotedblleft $\mathcal{%
\mathcal{B}}(x)=0$ only at $x=0$ \textquotedblright . Indeed, in the proof
of Theorem \ref{msb}, assume that the minimizing problem $\inf_{\{\mathcal{B}%
(x)=R\}}\mathcal{A}(x)$ is attained at $x_{R}\in X$ then because $\mathcal{A}
$ and $\mathcal{B}$ are differentiable there exists $(\lambda _{1},\lambda
_{2})$ a pair of Lagrange multipliers such that 
\begin{equation*}
\lambda _{1}\mathcal{A}^{\prime }(x_{R})+\lambda _{2}\mathcal{B}^{\prime
}(x_{R})=0.
\end{equation*}
Consequently, $\lambda _{1}$ and $\lambda _{2}$ are not both zero. In fact,
if $\lambda _{2}\neq 0$ and $\lambda _{1}=0$ then we get 
\begin{equation*}
\lambda _{2}(\mathcal{B}^{\prime }(x_{R}),x_{R})=0.
\end{equation*}
So, for instance, assume that the following condition obeys
\textquotedblleft there exists $\gamma >0$ such that 
\begin{equation*}
(\mathcal{B}^{\prime }(x),x)\geq \gamma \mathcal{B}(x)\hbox{ for all }x\in X.
\end{equation*}
In this case, particularly, taking $x=x_{R},$ it follows that $(\mathcal{B}%
^{\prime }(x_{R}),x_{R})=0$ implies $\mathcal{B}(x_{R})=0.$ This is a
contradiction because $x_{R}$ belongs in the set $\{\mathcal{B}(x)=R\}.$
\end{remark}


\subsection{ Properties on $\mathcal{A}$ and $\mathcal{B}$}

\begin{lemma}
\label{L3} ${(i)}$ $\mathcal{A}(z,w)$ is coercive on $X_0^{p(x),q(x)}(%
\Omega).$

$(ii)$ $\mathcal{B}$ is a weakly continuous functional, namely, $%
(z_{n},w_{n}){\rightharpoonup }(z,w)$ (weak convergence) implies $\mathcal{B}%
(z_{n},w_{n})\rightarrow \mathcal{B}(z,w)$.

$(iii)$ Let $(z,w)$ be in $X_0^{p(x),q(x)}(\Omega).$ Assume $\mathcal{B}%
^{\prime}(z,w)=0$ in $X^{-1,p^\prime(x),q^\prime(x)}(\Omega)$ then $\mathcal{%
B}(z,w)=0.$
\end{lemma}

\textbf{Proof.} \ignorespaces
(i)\,\, For any $(z,w)\in X_{0}^{p(x),q(x)}(\Omega )$ with $\left\Vert
z\right\Vert _{1,p(x)},\left\Vert w\right\Vert _{1,q(x)}>1$, using Lemma \ref%
{L1} we have 
\begin{equation*}
\begin{array}{l}
\int_{\Omega }\frac{1}{p(x)}\left\vert \nabla z\right\vert
^{p(x)}dx+\int_{\Omega }\frac{1}{q(x)}\left\vert \nabla w\right\vert
^{q(x)}dx \\ 
\geq \frac{1}{p^{+}}\int_{\Omega }\left\vert \nabla z\right\vert ^{p(x)}dx+%
\frac{1}{q^{+}}\int_{\Omega }\left\vert \nabla w\right\vert ^{q(x)}dx \\ 
\geq \min \{\frac{1}{p^{+}},\frac{1}{q^{+}}\}(\left\Vert z\right\Vert
_{1,p(x)}^{p^{-}}+\left\Vert w\right\Vert _{1,q(x)}^{q^{-}}) \\ 
\geq 2^{-\min \{p^{-},q^{-}\}}\min \{\frac{1}{p^{+}},\frac{1}{q^{+}}%
\}(\left\Vert z\right\Vert _{1,p(x)}+\left\Vert w\right\Vert
_{1,q(x)})^{\min \{p^{-},q^{-}\}}.%
\end{array}%
\end{equation*}
Since $\min \{p^{-},q^{-}\}>1$ (see \textrm{(H.3)} and \textrm{(H.4)}) the
above inequality implies that%
\begin{equation*}
\begin{array}{l}
\mathcal{A}(z,w)\rightarrow \infty \hbox{ \ as \ }\left\Vert
(z,w)\right\Vert \rightarrow \infty.%
\end{array}%
\end{equation*}

(ii)\,\, Let $(z_{n},w_{n})\rightharpoonup ({z,w)}$ in $X_{0}^{p(x),q(x)}(%
\Omega ).$ By the first part in \textrm{(H.4)} and (\ref{20}) the embeddings 
$W_{0}^{1,p(x)}\hookrightarrow L^{p(x)}(\Omega )$ and $W_{0}^{1,q(x)}%
\hookrightarrow L^{q(x)}(\Omega )$ are both compact, so we get 
\begin{equation}
(z_{n},w_{n})\rightarrow ({z,w)}\hbox{ in }L^{p(x)}(\Omega )\times
L^{q(x)}(\Omega ).  \label{21}
\end{equation}%
Using \textrm{(H.3)} and the definition of $\mathcal{B},$ we have%
\begin{equation*}
\begin{array}{rl}
\left\vert \mathcal{B}(z_{n},w_{n})-\mathcal{B}({z},{w})\right\vert \leq & 
\Vert c\Vert _{\infty }\left[ \int_{\Omega }|{z}|^{\alpha (x)+1}\left( |{w}%
|^{\beta (x)+1}-|{w_{n}}|^{\beta (x)+1}\right) dx\right. \\ 
& \left. +\int_{\Omega }|{w_{n}}|^{\alpha (x)+1}\left( |{z}|^{\alpha (x)+1}-|%
{z_{n}}|^{\alpha (x)+1}\right) dx\right] \\ 
\leq & 2^{\max (\alpha ^{+},\beta ^{+})}\Vert c\Vert _{\infty }\left[
\int_{\Omega }|{z}|^{\alpha (x)+1}\left\vert {w}-{w_{n}}\right\vert ^{\beta
(x)+1}dx\right. \\ 
& \left. +\int_{\Omega }|{\ w_{n}}|^{\alpha (x)+1}\left\vert {z}-{z_{n}}%
\right\vert ^{\alpha (x)+1}dx\right] .%
\end{array}%
\end{equation*}
By H\"{o}lder inequality one has 
\begin{equation*}
\begin{array}{l}
\int_{\Omega }|{z}|^{\alpha (x)+1}\left\vert {w}-{w_{n}}\right\vert ^{\beta
(x)+1}dx \\ 
\leq C_{\alpha ,\beta ,p,q}\left\Vert |{z}|^{\alpha (x)+1}\right\Vert _{L^{%
\frac{p(x)}{\alpha (x)+1}}(\Omega )}\left\Vert |w_{n}-{w}|^{\alpha
(x)+1}\right\Vert _{L^{\frac{q(x)}{\beta (x)+1}}(\Omega )}.%
\end{array}%
\end{equation*}
where $C_{\alpha ,\beta ,p,q}>0$ is a constant. Observe that%
\begin{equation*}
\begin{array}{l}
\left\Vert |{w}-{w_{n}}|^{\beta (x)+1}\right\Vert _{L^{\frac{q(x)}{\beta
(x)+1}}(\Omega )}^{q^{+}}\leq \int_{\Omega }(\left\vert {w}-{w_{n}}%
\right\vert ^{\beta (x)+1})^{\frac{q(x)}{\beta (x)+1}} dx=\rho _{q(.)}({w}-{%
w_{n}})%
\end{array}%
\end{equation*}
and%
\begin{equation*}
\begin{array}{l}
\rho _{q(.)}({w}-{w_{n}})\leq \left\Vert w-{w_{n}}|^{\beta (x)+1}\right\Vert
_{L^{{q(x)}}(\Omega )}^{q^{-}}.%
\end{array}%
\end{equation*}

Then it follows that 
\begin{equation*}
\begin{array}{l}
\left\Vert |{w}-{w_{n}}|^{\beta (x)+1}\right\Vert _{L^{\frac{q(x)}{\beta
(x)+1}}(\Omega )}\leq \rho _{q(.)}({w}-{w_{n}})^{1/q^{+}}\leq \left\Vert |{w}%
-{w_{n}}|^{\beta (x)+1}\right\Vert _{L^{{q(x)}}(\Omega )}^{q^{-}/q^{+}}.%
\end{array}%
\end{equation*}
Therefore, the strong convergence in (\ref{21}) ensures that 
\begin{equation*}
\begin{array}{l}
\left\Vert |{w}-{w_{n}}|^{\beta (x)+1}\right\Vert _{L^{\frac{q(x)}{\beta
(x)+1}}(\Omega )}\rightarrow 0\hbox{ \ as }n\rightarrow +\infty .%
\end{array}%
\end{equation*}
A quite similar argument provides 
\begin{equation*}
\begin{array}{l}
\left\Vert |{z-z_{n}}|^{\alpha (x)+1}\right\Vert _{L^{\frac{p(x)}{\alpha
(x)+1}}(\Omega )}\rightarrow 0\hbox{ \ as }n\rightarrow +\infty .%
\end{array}%
\end{equation*}

(iii)\,\, From (\ref{Bprime}), it is clear to notice that for any $(z,w)\in
X_0^{p(x),q(x)}(\Omega),$ by taking $\varphi=1/{p(x)}z$ and $\psi=1/{q(x)}w,$
the following identity holds 
\begin{equation*}
\mathcal{B}^{\prime}(z,w),(1/{p(x)}z,1/{q(x)}w)=\mathcal{B}(z,w).
\end{equation*}
Then the statement $(iii)$ follows. This conclude the proof of the Lemma. {\ 
\rule{0.2cm}{0.2cm}{\ }}


\subsection{A priori bound for $\mathcal{A}$}

\begin{lemma}
\label{L4} Let $R$ a fixed and strictly positive real. There exists a
constant $\mathcal{K}(R)>0$ depending on $R$ such that 
\begin{equation}
\begin{array}{l}
\mathcal{A}(z,w)\geq \mathcal{K}(R)>0,\,\,\, \forall (z,w)\in \mathcal{X}_R.%
\end{array}
\label{KR}
\end{equation}
\end{lemma}

\textbf{Proof.} \ignorespaces
First, observe from Lemma \ref{L2} that if $\Vert \nabla z\Vert
_{L^{p(x)}(\Omega )}<1,$ we have 
\begin{equation*}
\begin{array}{l}
\left\Vert \frac{z}{C_{N,p}}\right\Vert _{L^{p(x)}(\Omega )}<1.%
\end{array}%
\end{equation*}
Then if follows that 
\begin{equation}
\begin{array}{l}
\rho _{p(x)}(\frac{z}{C_{N,p}})\leq \left\Vert \frac{z}{C_{N,p}}\right\Vert
_{L^{p(x)}(\Omega )}^{p^{-}},%
\end{array}
\label{2}
\end{equation}%
which combined with Lemma \ref{L2} leads to%
\begin{equation*}
\begin{array}{l}
\int_{\Omega }\frac{|z|^{p(x)}}{C_{N,p}^{p(x)}}dx\leq \Vert \nabla z\Vert
_{L^{p(x)}(\Omega )}^{p^{-}}.%
\end{array}%
\end{equation*}
Hence it holds 
\begin{equation}
\begin{array}{l}
\int_{\Omega }|z|^{p(x)}dx\leq K_{N,p}\Vert \nabla z\Vert _{L^{p(x)}(\Omega
)}^{p^{-}}\leq K_{N,p}\Vert \nabla z\Vert _{L^{p(x)}(\Omega )}^{p^{-}/p^{+}},%
\end{array}
\label{3}
\end{equation}%
where 
\begin{equation*}
K_{N,p}=\left\{ 
\begin{array}{rl}
C_{N,p}^{p^{+}} & \hbox{ if \ }C_{N,p}>1 \\ 
C_{N,p}^{p^{-}} & \hbox{ if \ }C_{N,p}<1.%
\end{array}%
\right.
\end{equation*}
A quite similar argument shows that%
\begin{equation}
\begin{array}{l}
\int_{\Omega }|w|^{q(x)}dx\leq K_{N,q}\Vert \nabla w\Vert _{L^{q(x)}(\Omega
)}^{q^{-}/q^{+}},%
\end{array}
\label{4}
\end{equation}%
where 
\begin{equation*}
K_{N,q}=\left\{ 
\begin{array}{rl}
C_{N,q}^{q^{+}} & \hbox{ if \ }C_{N,q}>1 \\ 
C_{N,q}^{q^{-}} & \hbox{ if \ }C_{N,q}<1.%
\end{array}%
\right.
\end{equation*}
For every $(z,w)\in X_{0}^{p(x),q(x)}(\Omega ),$ Young inequality and 
\textrm{(H.3)} imply 
\begin{equation}
\begin{array}{rl}
\int_{\Omega }c(x)|z|^{\alpha (x)+1}|w|^{\beta (x)+1}dx & \leq \Vert c\Vert
_{\infty }\int_{\Omega }\left[ \frac{\alpha (x)+1}{p(x)}|z|^{p(x)}+\frac{%
\beta (x)+1}{q(x)}|w|^{q(x)}\right] dx \\ 
& \leq \Vert c\Vert _{\infty }(\int_{\Omega }|z|^{p(x)}dx+\int_{\Omega
}|w|^{q(x)}dx).%
\end{array}
\label{6}
\end{equation}%
Assume that $(z,w)\in \mathcal{X}_{R}$ is such that 
\begin{equation}
\begin{array}{l}
\max \left( \Vert \nabla z\Vert _{L^{p(\cdot )}(\Omega )},\Vert \nabla
w\Vert _{L^{q(\cdot )}(\Omega )}\right) <1.%
\end{array}
\label{5}
\end{equation}%
Bearing in mind \textrm{(H.3)}, \textrm{(H.4)} and ${(i)}$ of Lemma \ref{L1}%
, we have 
\begin{equation}
\begin{array}{l}
\max \left\{ \int_{\Omega }\frac{1}{p(x)}|\nabla z|^{p(x)}dx,\hbox{ }%
\int_{\Omega }\frac{1}{q(x)}|\nabla w|^{q(x)}dx\right\} <1.%
\end{array}
\label{7}
\end{equation}%
Then, from (\ref{3})-(\ref{7}), it follows that%
\begin{equation}
\begin{array}{l}
R\leq {K}_{1}\left( \int_{\Omega }\frac{1}{p(x)}|\nabla z|^{p(x)}dx\right)
^{p^{-}/p^{+}}+{K}_{2}\left( \int_{\Omega }\frac{1}{q(x)}|\nabla
w|^{q(x)}dx\right) ^{q^{-}/q^{+}}%
\end{array}
\label{8}
\end{equation}%
From the hypothesis \textrm{(H.4)} on $p^{-},$ $p^{+},$ $q^{-}$ and $q^{+},$
it follows that 
\begin{equation}
\begin{array}{l}
R^{\frac{p^{+}q^{+}}{p^{-}q^{-}}}\leq 2^{\frac{p^{+}q^{+}}{p^{-}q^{-}}-1}%
\left[ {K}_{1}^{\frac{p^{+}q^{+}}{p^{-}q^{-}}}\left( \int_{\Omega }\frac{1}{%
p(x)}|\nabla z|^{p(x)}dx\right) ^{q^{+}/q^{-}}\right. \\ 
\left. \qquad \qquad \qquad \qquad \qquad +{K}_{2}^{p^{+}q^{+}/p^{-}q^{-}}%
\left( \int_{\Omega }\frac{1}{q(x)}|\nabla w|^{q(x)}dx\right) ^{p^{+}/p^{-}}%
\right] .%
\end{array}
\label{9}
\end{equation}%
Or again 
\begin{equation}
R^{\frac{p^{+}q^{+}}{p^{-}q^{-}}}\leq (2{K}_{3})^{\frac{p^{+}q^{+}}{%
p^{-}q^{-}}}\left[ \int_{\Omega }\frac{1}{p(x)}|\nabla
z|^{p(x)}dx+\int_{\Omega }\frac{1}{q(x)}|\nabla w|^{q(x)}dx\right]
\label{10}
\end{equation}%
where%
\begin{equation*}
\begin{array}{l}
{K}_{1}=K_{N,p}\left( p^{+}\right) ^{p^{-}/p^{+}}\Vert c\Vert _{\infty },{\
\ }{K}_{2}=K_{N,q}\left( q^{+}\right) ^{q^{-}/q^{+}}\Vert c\Vert _{\infty }%
\end{array}%
\end{equation*}
and ${K}_{3}={K}_{1}+{K}_{2}.$ Thus, from (\ref{10}), we conclude that 
\begin{equation}
\begin{array}{l}
\mathcal{A}(z,w)\geq \left( \frac{R}{2{K}_{3}}\right) ^{\frac{q^{+}p^{+}}{%
q^{-}p^{-}}}.%
\end{array}
\label{11}
\end{equation}%
Now, we deal with the case when $(z,w)\in \mathcal{X}_{R}$ is such that 
\begin{equation*}
\begin{array}{l}
\max \left( \Vert \nabla z\Vert _{L^{p(\cdot )}(\Omega )},\Vert \nabla
w\Vert _{L^{q(\cdot )}(\Omega )}\right) \geq 1.%
\end{array}%
\end{equation*}
This implies that%
\begin{equation*}
\begin{array}{l}
\max \left( \int_{\Omega }|\nabla z|^{p(x)}dx,\int_{\Omega }|\nabla
w|^{q(x)}dx\right) \geq 1.%
\end{array}%
\end{equation*}
If $\int_{\Omega }|\nabla z|^{p(x)}dx\geq 1$ we have%
\begin{equation*}
\begin{array}{l}
p^{+}\int_{\Omega }\frac{1}{p(x)}|\nabla z|^{p(x)}dx\geq \int_{\Omega
}|\nabla z|^{p(x)}dx\geq 1.%
\end{array}%
\end{equation*}
which in turn yields 
\begin{equation}
\begin{array}{l}
\mathcal{A}(z,w)=\int_{\Omega }\frac{1}{p(x)}|\nabla
z|^{p(x)}dx+\int_{\Omega }\frac{1}{q(x)}|\nabla w|^{q(x)}dx>\frac{1}{p^{+}}.%
\end{array}
\label{12}
\end{equation}%
Now for $\int_{\Omega }|\nabla w|^{q(x)}dx\geq 1$ a quite similar argument
provides 
\begin{equation}
\begin{array}{l}
\mathcal{A}(z,w)>\frac{1}{q^{+}}.%
\end{array}
\label{13}
\end{equation}%
We notice that if $\max \left( \Vert \nabla z\Vert _{L^{p(\cdot )}(\Omega
)},\Vert \nabla w\Vert _{L^{q(\cdot )}(\Omega )}\right) \geq 1,$ from (\ref%
{12}) and (\ref{13}), it is clearly that 
\begin{equation}
\begin{array}{l}
\mathcal{A}(z,w)>\max (\frac{1}{p^{+}},\frac{1}{q^{+}}).%
\end{array}
\label{14}
\end{equation}

Thus, according to (\ref{11}) and (\ref{14}), for all $(z,w)\in \mathcal{X}%
_{R},$ one has%
\begin{equation}
\begin{array}{l}
\mathcal{A}(z,w)\geq \max \{(\frac{R}{2{K}_{3}})^{\frac{q^{+}p^{+}}{%
q^{-}p^{-}}},\frac{1}{p^{+}},\frac{1}{q^{+}}\}>0.%
\end{array}
\label{15}
\end{equation}%
Consequently, there exists a constant $\mathcal{K}(R)>0$ depending on $R$
such that (\ref{KR}) holds. {\ \rule{0.2cm}{0.2cm}{\ }}

\subsection{Proof of (\protect\ref{infAR})}

\label{SS2} We begin by the proposition.

\begin{proposition}
\label{T1}Assume that \textrm{(H.3)} holds. Then, for $R>0,$

\begin{description}
\item[\textrm{(i)}] $0<\frac{\lambda ^{\ast }_R}{(\alpha^++\beta^++2)}< {%
\lambda _{\ast }}_R<\lambda ^{\ast }_R.$

\item[\textrm{(ii)}] Any $\lambda <{\lambda _{\ast }}_R$ is not an
eigenvalue of problem (\ref{p}).

\item[\textrm{(iii)}] There exists $(\hat{u}_{R},\hat{v}_{R})\in \mathcal{X}%
_{R}$ such that $\lambda _{R}^{\ast }$ is a corresponding eigenvalue for the
system (\ref{p}).
\end{description}
\end{proposition}

\textbf{Proof.} \ignorespaces
\textbf{(i)}. First let us show that $0<\frac{\lambda _{R}^{\ast }}{(\alpha
^{+}+\beta ^{+}+2)}\leq {\lambda _{\ast }}_{R}\leq \lambda _{R}^{\ast }.$
Obviously, for all $(z,w)\in \mathcal{X}_{R},$ we have 
\begin{equation*}
\begin{array}{l}
\frac{\mathcal{A}(z,w)}{(\alpha ^{+}+\beta ^{+}+2)R}\leq \frac{\mathcal{A}%
(z,w)}{\int_{\Omega }c(x)(\alpha (x)+\beta (x)+2)|z|^{\alpha
(x)+1}|w|^{\beta (x)+1}dx}\leq \frac{\mathcal{A}(z,w)}{R}.%
\end{array}%
\end{equation*}
from (\ref{c1}) and (\ref{c2}), it derives that $\frac{\lambda _{R}^{\ast }}{%
(\alpha ^{+}+\beta ^{+}+2)}<{\lambda _{\ast }}_{R}<\lambda _{R}^{\ast }.$
Now suppose that ${\lambda _{\ast }}_{R}=0$. Then $\lambda _{R}^{\ast }=0$
and by virtue of Lemma \ref{L4}\ and Remark \ref{R1} this is a
contradiction. Hence ${\lambda _{\ast }}_{R}>0$. \medskip

\textbf{(ii)}. Next we show that $\lambda $ cannot be an eigenvalue for $%
\lambda <\lambda _{\ast }$. Indeed, suppose by contradiction that $\lambda $
is an eigenvalue of problem (\ref{p}). Then there exists $(u,v)\in
X_{0}^{p(x),q(x)}(\Omega )-\{(0,0)\}$ such that%
\begin{equation}
\begin{array}{l}
\int_{\Omega }\left\vert \nabla u\right\vert ^{p(x)} dx=\lambda \int_{\Omega
}c(x)(\alpha (x)+1)|u|^{\alpha (x)+1}|v|^{\beta (x)+1} \\ 
\int_{\Omega }\left\vert \nabla v\right\vert ^{q(x)} dx=\lambda \int_{\Omega
}c(x)(\beta (x)+1)|u|^{\alpha (x)+1}|v|^{\beta (x)+1}.%
\end{array}
\label{35}
\end{equation}%
On the basis of \textrm{(H.3)}, \textrm{(H.4)}, (\ref{c2}) and (\ref{35}),
we get%
\begin{equation*}
\begin{array}{l}
\lambda _{\ast }\int_{\Omega }c(x)(\alpha (x)+\beta (x)+2)|u|^{\alpha
(x)+1}|v|^{\beta (x)+1}dx \\ 
\leq \int_{\Omega }\left( \frac{1}{p(x)}\left\vert \nabla u\right\vert
^{p(x)}+\frac{1}{q(x)}\left\vert \nabla v\right\vert ^{q(x)}\right) dx \\ 
\leq \int_{\Omega }\left\vert \nabla u\right\vert ^{p(x)}dx+\int_{\Omega
}\left\vert \nabla v\right\vert ^{q(x)}dx \\ 
=\lambda \int_{\Omega }c(x)(\alpha (x)+\beta (x)+2)|u|^{\alpha
(x)+1}|v|^{\beta (x)+1}dx \\ 
<\lambda _{\ast }\int_{\Omega }c(x)(\alpha (x)+\beta (x)+2)|u|^{\alpha
(x)+1}|v|^{\beta (x)+1}dx,%
\end{array}%
\end{equation*}
which is not possible and the conclusion follows. \medskip

\textbf{(iii)}. Now, we claim that the infimum in (\ref{infAR}) is achieved
at an element of $\mathcal{X}_{R}.$ Indeed, thanks to the lemma \ref{L3}, $%
\mathcal{B}$ is weakly continuous on $X_{0}^{p(x),q(x)}(\Omega ),$ then the
nonempty set $\mathcal{X}_{R}$ is weakly closed. So, since $\mathcal{A}$ is
weakly lower semicontinuous, we conclude that there exists an element of $%
\mathcal{X}_{R}$ which we denote $\left( {\hat{u}},{\hat{v}_{R}}\right) $
such that (\ref{infAR}) is feasible. Since $({\hat{u}_{R}},{\hat{v}_{R}}%
)\neq 0,$ we also have $\mathcal{B}^{\prime }\left( {\hat{u}_{R}},{\hat{v}%
_{R}}\right) \neq 0$ otherwise it implies $\mathcal{B}\left( {\hat{u}_{R}},{%
\hat{v}_{R}}\right) =0$ and which contradicts $\left( {\hat{u}},{\hat{v}_{R}}%
\right) \in \mathcal{X}_{R}.$ So, owing to Lagrange multiplier method (see
e.g. \cite[Theorem 6.3.2, p. 325]{msb} or \cite[Theorem 6.3.2, p. 402]{DM}),
there exists ${\lambda _{R}}\in 
\mathbb{R}
$ such that%
\begin{equation}
\mathcal{A}^{\prime }({\hat{u}_{R}},{\hat{v}_{R}})\cdot (\varphi ,\psi
)=\lambda _{R}\mathcal{B}^{\prime }({\hat{u}_{R}},{\hat{v}_{R}})\cdot
(\varphi ,\psi ),\,\,\,\forall (\varphi ,\psi )\in X_{0}^{p(x),q(x)}(\Omega )
\label{56}
\end{equation}%
where $\mathcal{A}^{\prime }$ and $\mathcal{B}^{\prime }$ are defined as in (%
\ref{Aprime}) and (\ref{Bprime}) respectively.\medskip

In the sequel, we show that $\lambda _{R}$ is equal to $\lambda _{R}^{\ast
}. $ To this end, let us denote by $\Omega ^{+}$ and $\Omega ^{-}$ the sets
defined as follows 
\begin{equation*}
\begin{array}{l}
\Omega ^{+}=\{x\in \Omega ;\,\,|\nabla {\hat{u}_{R}}|^{p(x)}-{\lambda _{R}}%
(\alpha (x)+1)c(x)|{\hat{u}_{R}}|^{\alpha (x)+1}|{\hat{v}_{R}}|^{\beta
(x)+1}\geq 0\}%
\end{array}%
\end{equation*}%
and%
\begin{equation*}
\begin{array}{l}
\Omega _{-}=\{x\in \Omega ;\,\,|\nabla {\hat{u}_{R}}|^{p(x)}-{\lambda _{R}}%
(\alpha (x)+1)c(x)|{\hat{u}_{R}}|^{\alpha (x)+1}|{\hat{v}_{R}}|^{\beta
(x)+1}<0\}.%
\end{array}%
\end{equation*}%
By taking $\varphi ={\hat{u}_{R}}$ $1_{\Omega ^{+}}$ and $\psi =0$ in (\ref%
{56}) one has 
\begin{equation}
\begin{array}{l}
\int_{\Omega ^{+}}\left( |\nabla {\hat{u}_{R}}|^{p(x)}-{\lambda _{R}}%
c(x)(\alpha (x)+1)|{\hat{u}_{R}}|^{\alpha (x)+1}|{\hat{v}_{R}}|^{\beta
(x)+1}\right) dx=0%
\end{array}
\label{58*}
\end{equation}%
and likewise, by choosing $\varphi ={\hat{u}_{R}}$ $1_{\Omega ^{-}}$ and $%
\psi =0$ in (\ref{56}) we get 
\begin{equation}
\begin{array}{l}
\int_{\Omega ^{-}}\left( |\nabla {\hat{u}_{R}}|^{p(x)}-{\lambda _{R}}%
c(x)(\alpha (x)+1)|{\hat{u}_{R}}|^{\alpha (x)+1}|{\hat{v}_{R}}|^{\beta
(x)+1}\right) dx=0.%
\end{array}
\label{58}
\end{equation}%
We claim that%
\begin{equation}
\begin{array}{l}
\int_{\Omega }\frac{1}{p(x)}|\nabla {\hat{u}_{R}}|^{p(x)}dx={\lambda _{R}}%
\int_{\Omega }c(x)\frac{\alpha (x)+1}{p(x)}|{\hat{u}_{R}}|^{\alpha (x)+1}|{%
\hat{v}}|^{\beta (x)+1}dx.%
\end{array}
\label{59}
\end{equation}%
Indeed, on account of \textrm{(H.4)}, (\ref{58*}) and (\ref{58}) we have%
\begin{equation*}
\begin{array}{l}
\left\vert \int_{\Omega }\frac{|\nabla {\hat{u}_{R}}|^{p(x)}}{p(x)}dx-{%
\lambda _{R}}\int_{\Omega }\frac{\alpha (x)+1}{p(x)}c(x)|{\hat{u}_{R}}%
|^{\alpha (x)+1}|{\hat{v}}|^{\beta (x)+1}dx\right\vert \\ 
\leq \int_{\Omega }p(x)\left\vert \frac{|\nabla {\hat{u}_{R}}|^{p(x)}}{p(x)}-%
{\lambda _{R}}\frac{\alpha (x)+1}{p(x)}c(x)|{\hat{u}_{R}}|^{\alpha (x)+1}|{%
\hat{v}_{R}}|^{\beta (x)+1}\right\vert dx \\ 
=\int_{\Omega }\left\vert |\nabla {\hat{u}_{R}}|^{p(x)}-{\lambda _{R}}%
(\alpha (x)+1)c(x)|{\hat{u}_{R}}|^{\alpha (x)+1}|{\hat{v}_{R}}|^{\beta
(x)+1}\right\vert dx \\ 
\leq \int_{\Omega ^{+}}\left( |\nabla {\hat{u}_{R}}|^{p(x)}-{\lambda _{R}}%
(\alpha (x)+1)c(x)|{\hat{u}_{R}}|^{\alpha (x)+1}|{\hat{v}_{R}}|^{\beta
(x)+1}\right) dx \\ 
-\int_{\Omega _{-}}\left( |\nabla {\hat{u}_{R}}|^{p(x)}-{\lambda _{R}}%
(\alpha (x)+1)c(x)|{\hat{u}_{R}}|^{\alpha (x)+1}|{\hat{v}_{R}}|^{\beta
(x)+1}\right) dx=0,%
\end{array}%
\end{equation*}%
showing that (\ref{59}) holds. In the same manner we can prove that%
\begin{equation}
\begin{array}{l}
\int_{\Omega }\frac{1}{q(x)}|\nabla {\hat{v}_{R}}|^{q(x)}dx={\lambda _{R}}%
\int_{\Omega }c(x)\frac{\beta (x)+1}{q(x)}|{\hat{u}_{R}}|^{\alpha (x)+1}|{%
\hat{v}_{R}}|^{\beta (x)+1}dx.%
\end{array}
\label{59*}
\end{equation}%
Adding together (\ref{59}) and (\ref{59*}), on account of \textrm{(H.3)} and
(\ref{13}), we achieve that%
\begin{equation*}
\mathcal{A}({\hat{u}_{R}},{\hat{v}_{R}})=R{\lambda _{R}}.
\end{equation*}%
Then, bearing in mind (\ref{14}) it turns out that ${\lambda _{R}}=\lambda
_{R}^{\ast }$, showing that $\lambda _{R}^{\ast }$ is at least one
eigenvalue of (\ref{p}).

Then, combining this last point with the characterization (\ref{56}), we get%
\begin{equation*}
\mathcal{A}^{\prime }({\hat{u}_{R}},{\hat{v}_{R}})\cdot (\varphi ,0)=\lambda
_{R}^{\ast }\mathcal{B}^{\prime }({\hat{u}_{R}},{\hat{v}_{R}})\cdot (\varphi
,0),\,\,\,\forall \varphi \in W_{0}^{1,q(x)}(\Omega )
\end{equation*}
and 
\begin{equation*}
\mathcal{A}^{\prime }({\hat{u}_{R}},{\hat{v}_{R}})\cdot (0,\psi )=\lambda
_{R}^{\ast }\mathcal{B}^{\prime }({\hat{u}_{R}},{\hat{v}_{R}})\cdot (0,\psi
),\,\,\,\forall \psi \in W_{0}^{1,q(x)}(\Omega ).
\end{equation*}
On other words, it means that $(({\hat{u}_{R}},{\hat{v}_{R}}),\lambda
_{R}^{\ast })$ is a solution of the system (\ref{p}). {\ \rule{0.2cm}{0.2cm}{%
\ }}

\subsection{Proof of Theorem \protect\ref{T0}}

Employing again the statement of Lemma \ref{L3}, we can apply the theorem %
\ref{msb} due to \cite{msb}. Then the system (\ref{p}) has a one-parameter
family of nontrivial solutions $(({\hat{u}_{R}},{\hat{v}_{R}}),\lambda _{R})$
for all $R>0.$ Moreover, from \textbf{(iii)} of Proposition \ref{T1}, $%
\lambda _{R}=\lambda _{R}^{\ast }.$

It remains to prove that $\lambda _{p(x),q(x)}^{\ast }=\inf_{R>0}\lambda
_{R}^{\ast }>0$. From (\ref{6}) and for $(z,w)\in X_{0}^{1,p(x),q(x)}(\Omega
)\setminus \{0\}$, one has%
\begin{equation}
\begin{array}{l}
\frac{1}{\Vert c\Vert _{\infty }}\cdot \frac{\int_{\Omega }\frac{1}{p(x)}%
\left\vert \nabla z\right\vert ^{p(x)}dx+\int_{\Omega }\frac{1}{q(x)}%
\left\vert \nabla w\right\vert ^{q(x)}dx}{\int_{\Omega }\left\vert
z\right\vert ^{p(x)}dx+\int_{\Omega }\left\vert w\right\vert ^{q(x)}dx}\leq 
\frac{\int_{\Omega }\frac{1}{p(x)}\left\vert \nabla z\right\vert
^{p(x)}dx+\int_{\Omega }\frac{1}{q(x)}\left\vert \nabla w\right\vert
^{q(x)}dx}{\int_{\Omega }c(x)|z|^{\alpha (x)+1}|w|^{\beta (x)+1}dx}.%
\end{array}
\label{e}
\end{equation}%
Recalling that under assumption \textrm{(a.1)} or \textrm{(a.2)}, the
authors in \cite{fanandall} proved that the first eignevalues 
\begin{equation}
\left\{ 
\begin{array}{l}
\lambda _{p(x)}^{\ast }=\inf_{z\in W_{0}^{1,p(x)}(\Omega )\setminus \{0\}}%
\frac{\int_{\Omega }\left\vert \nabla z\right\vert ^{p(x)}dx}{\int_{\Omega
}\left\vert z\right\vert ^{p(x)}dx} \\ 
\lambda _{q(x)}^{\ast }=\inf_{z\in W_{0}^{1,q(x)}(\Omega )\setminus \{0\}}%
\frac{\int_{\Omega }\left\vert \nabla z\right\vert ^{q(x)}dx}{\int_{\Omega
}\left\vert z\right\vert ^{q(x)}dx},%
\end{array}%
\right.  \label{l*}
\end{equation}%
are strictly positive. Hence, combining with (\ref{e}) it follows that%
\begin{equation*}
\begin{array}{l}
\min \left\{ \frac{\lambda _{p(x)}^{\ast }}{p^{+}\Vert c\Vert _{\infty }},%
\frac{\lambda _{q(x)}^{\ast }}{q^{+}\Vert c\Vert _{\infty }}\right\} =\min
\left\{ \frac{\lambda _{p(x)}^{\ast }}{p^{+}\Vert c\Vert _{\infty }},\frac{%
\lambda _{q(x)}^{\ast }}{q^{+}\Vert c\Vert _{\infty }}\right\} \cdot \frac{%
\int_{\Omega }\left\vert z\right\vert ^{p(x)}dx+\int_{\Omega }\left\vert
w\right\vert ^{q(x)}dx}{\int_{\Omega }\left\vert z\right\vert
^{p(x)}dx+\int_{\Omega }\left\vert w\right\vert ^{q(x)}dx} \\ 
\\ 
\leq \frac{\frac{\lambda _{p(x)}^{\ast }}{p^{+}\Vert c\Vert _{\infty }}%
\int_{\Omega }\left\vert z\right\vert ^{p(x)}dx+\frac{\lambda _{q(x)}^{\ast }%
}{q^{+}\Vert c\Vert _{\infty }}\int_{\Omega }\left\vert w\right\vert
^{q(x)}dx}{\int_{\Omega }\left\vert z\right\vert ^{p(x)}dx+\int_{\Omega
}\left\vert w\right\vert ^{q(x)}dx}\leq \frac{\int_{\Omega }\frac{1}{p(x)}%
\left\vert \nabla z\right\vert ^{p(x)}dx+\int_{\Omega }\frac{1}{q(x)}%
\left\vert \nabla w\right\vert ^{q(x)}dx}{\int_{\Omega }c(x)|z|^{\alpha
(x)+1}|w|^{\beta (x)+1}dx}.%
\end{array}%
\end{equation*}
Then 
\begin{equation}
\begin{array}{l}
0<\min \left\{ \frac{\lambda _{p(x)}^{\ast }}{p^{+}\Vert c\Vert _{\infty }},%
\frac{\lambda _{q(x)}^{\ast }}{q^{+}\Vert c\Vert _{\infty }}\right\} \\ 
\leq \inf_{(z,w)\in X_{0}^{1,p(x),q(x)}(\Omega )\setminus \{0\}}\frac{%
\int_{\Omega }\frac{1}{p(x)}\left\vert \nabla z\right\vert
^{p(x)}dx+\int_{\Omega }\frac{1}{q(x)}\left\vert \nabla w\right\vert
^{q(x)}dx}{\int_{\Omega }c(x)|z|^{\alpha (x)+1}|w|^{\beta (x)+1}dx}.%
\end{array}
\label{f}
\end{equation}%
Another hand, since 
\begin{equation*}
\begin{array}{l}
\bigcup_{R>0}\mathcal{X}_{R}\subset \left\{ (z,w)\in
X_{0}^{1,p(x),q(x)}(\Omega )\setminus \{0\}\right\} ,%
\end{array}%
\end{equation*}
one gets%
\begin{equation}
\begin{array}{l}
\inf_{(z,w)\in X_{0}^{1,p(x),q(x)}(\Omega )\setminus \{0\}}\frac{%
\int_{\Omega }\frac{1}{p(x)}\left\vert \nabla z\right\vert
^{p(x)}dx+\int_{\Omega }\frac{1}{q(x)}\left\vert \nabla w\right\vert
^{q(x)}dx}{\int_{\Omega }c(x)|z|^{\alpha (x)+1}|w|^{\beta (x)+1}dx} \\ 
\\ 
\leq \inf_{\{\mathcal{B}(z,w)=R\}}\frac{\int_{\Omega }\frac{1}{p(x)}%
\left\vert \nabla z\right\vert ^{p(x)}dx+\int_{\Omega }\frac{1}{q(x)}%
\left\vert \nabla w\right\vert ^{q(x)}dx}{\int_{\Omega }c(x)|z|^{\alpha
(x)+1}|w|^{\beta (x)+1}dx}.%
\end{array}
\label{g}
\end{equation}%
Thus, gathering (\ref{f}) and (\ref{g}) together we infer that 
\begin{equation*}
\begin{array}{l}
0<\min \left\{ \frac{\lambda _{p(x)}^{\ast }}{p^{+}\Vert c\Vert _{\infty }},%
\frac{\lambda _{q(x)}^{\ast }}{q^{+}\Vert c\Vert _{\infty }}\right\} \leq
\lambda _{p(x),q(x)}^{\ast }\leq \inf_{R>0}\lambda _{R}^{\ast }.%
\end{array}%
\end{equation*}
Next, let us prove that $\lambda _{p(x),q(x)}^{\ast }\geq \inf_{R>0}\lambda
_{R}^{\ast }$. To this end, let a constant $\varepsilon >0$, there is $%
R_{\varepsilon }>0$ such that $\lambda _{R_{\varepsilon }}^{\ast
}<\inf_{R>0}\lambda _{R}^{\ast }+\varepsilon .$ This implies that%
\begin{equation}
\lambda _{R_{\varepsilon }}^{\ast }<\lambda _{R}^{\ast }+\varepsilon 
\hbox{
\ for all }R>0\hbox{\ and }\varepsilon >0.  \label{63}
\end{equation}%
Now, let $(z,w)\in X_{0}^{1,p(x),q(x)}(\Omega )\setminus \{0\}$ such that $%
\mathcal{B}(z,w)>0$ and assume that $R_{(z,w)}=\mathcal{B}(z,w)$. According
to \textrm{(iii)} in Propostion \ref{T1}, the constant 
\begin{equation*}
\begin{array}{l}
\lambda _{R_{(z,w)}}^{\ast }=\inf_{\{\mathcal{B}(z,w)=R_{(z,w)}\}}\frac{%
\int_{\Omega }\frac{1}{p(x)}\left\vert \nabla z\right\vert
^{p(x)}dx+\int_{\Omega }\frac{1}{q(x)}\left\vert \nabla w\right\vert
^{q(x)}dx}{\int_{\Omega }c(x)|z|^{\alpha (x)+1}|w|^{\beta (x)+1}dx}%
\end{array}%
\end{equation*}
exists and then 
\begin{equation*}
\begin{array}{l}
\lambda _{R_{(z,w)}}^{\ast }\leq \frac{\int_{\Omega }\frac{1}{p(x)}%
\left\vert \nabla z\right\vert ^{p(x)}dx+\int_{\Omega }\frac{1}{q(x)}%
\left\vert \nabla w\right\vert ^{q(x)}dx}{\int_{\Omega }c(x)|z|^{\alpha
(x)+1}|w|^{\beta (x)+1}dx}.%
\end{array}%
\end{equation*}
At this point, combining with (\ref{63}) yields%
\begin{equation*}
\begin{array}{l}
\lambda _{R_{\varepsilon }}^{\ast }<\lambda _{R_{(z,w)}}^{\ast }+\varepsilon
\leq \frac{\mathcal{A}(z,w)}{\mathcal{B}(z,w)}+\varepsilon \hbox{ \ for all }%
\varepsilon >0,%
\end{array}%
\end{equation*}
which, it turn, leads to%
\begin{equation*}
\begin{array}{l}
\lambda _{R_{\varepsilon }}^{\ast }<\lambda _{R_{(z,w)}}^{\ast }+\varepsilon
\leq \inf_{(z,w)\in X_{0}^{1,p(x),q(x)}(\Omega )\setminus \{0\}}\frac{%
\mathcal{A}(z,w)}{\mathcal{B}(z,w)}+\varepsilon \hbox{ \ for all }%
\varepsilon >0.%
\end{array}%
\end{equation*}
This is equivalent to $\lambda _{R_{\varepsilon }}^{\ast }\leq \lambda
_{p(x),q(x)}^{\ast }+\varepsilon $. Consequently, 
\begin{equation*}
\begin{array}{l}
\inf_{R>0}\lambda _{R}^{\ast }\leq \lambda _{R_{\varepsilon }}^{\ast }\leq
\lambda _{p(x),q(x)}^{\ast }+\varepsilon \leq \inf_{R>0}\lambda _{R}^{\ast
}+\varepsilon \hbox{ \ for all }\varepsilon >0.%
\end{array}%
\end{equation*}
Finally, passing to the limit as $\varepsilon \rightarrow 0$ implies that $%
\lambda _{p(x),q(x)}^{\ast }=\inf_{R>0}\lambda _{R}^{\ast }$. This ends the
proof of Theorem \ref{T0}.

\section{Proof of Theorem \protect\ref{T2}}

\label{S4}

Let $(\hat{u}_{R},\hat{v}_{R})\in X_{0}^{p(x),q(x)}(\Omega )$ be a solution
of problem (\ref{p}) corresponding to the positive infimum eigenvalue $%
\lambda _{R}^{\ast }$ and let $d>0$ be a cosntant such that%
\begin{equation}
\begin{array}{l}
d=\frac{\hat{d}}{\max \left\{ p^{+},q^{+}\right\} },%
\end{array}
\label{53}
\end{equation}%
where%
\begin{equation}
\begin{array}{l}
1<\max \left\{ p^{+},q^{+}\right\} <\hat{d}\leq \max \left\{
p^{+},q^{+}\right\} \cdot \min \left\{ \frac{\pi _{p}^{-}}{p^{-}},\frac{\pi
_{p}^{+}}{p^{+}},\frac{\pi _{q}^{-}}{q^{-}},\frac{\pi _{q}^{+}}{q^{+}}%
\right\}%
\end{array}
\label{50}
\end{equation}%
and 
\begin{equation}
\begin{array}{l}
\pi _{p}(x)=\frac{Np(x)}{N-p(x)}, \, \pi _{p}^{-}=\inf_{x\in \Omega }\pi
_{p}(x)\hbox{ \ and \ }\pi _{p}^{+}=\sup_{x\in \Omega }\pi _{p}(x).%
\end{array}%
\end{equation}%
In this section, the goal consists in proving that $({\hat{u}_{R}},{\hat{v}%
_{R}})$ is bounded in $\Omega .$ Notice that from the above section, we have 
\begin{equation}
\left\{ 
\begin{array}{l}
\int_{\Omega }\left\vert \nabla {\hat{u}_{R}}\right\vert ^{p(x)-2}\nabla {%
\hat{u}_{R}}\nabla {\varphi } dx=\lambda _{R}^{\ast }\int_{\Omega
}c(x)(\alpha (x)+1){\hat{u}_{R}}|{\hat{u}_{R}}|^{\alpha (x)-1}|{\hat{v}_{R}}%
|^{\beta (x)+1}{\varphi } dx \\ 
\int_{\Omega }\left\vert \nabla {\hat{v}_{R}}\right\vert ^{q(x)-2}\nabla {%
\hat{v}_{R}}\nabla {\psi } dx=\lambda _{R}^{\ast }\int_{\Omega }c(x)(\beta
(x)+1)|{\hat{u}_{R}}|^{\alpha (x)+1}{\hat{v}_{R}}|{\hat{v}_{R}}|^{\beta
(x)-1}\psi dx.%
\end{array}%
\right.  \label{1}
\end{equation}

\begin{remark}
From the density of $C_{c}^{\infty }(\Omega )$ in $W_{0}^{1,p(x)}(\Omega )$
and through the embeddings $C_{c}^{\infty }(\Omega )\subset C^{1}(\overline{%
\Omega }),$ $C^{1}(\overline{\Omega })\subset W_{0}^{1,p^{+}}(\Omega )$ and $%
W_{0}^{1,p^{+}}(\Omega )\subset W_{0}^{1,p(x)}(\Omega )$ (since $p(x)\leq
p^{+}$ in $\Omega $), we may assume that ${\hat{u}_R}\in C^{1}(\overline{%
\Omega })$ (see, e.g., \cite{fz}). The same argument enable us to assume
that ${\hat{v}_R}\in C^{1}(\overline{\Omega })$.
\end{remark}

For a better reading, we divide the proof of Theorem \ref{T2} in several
lemmas.

\begin{lemma}
\label{estimk} Assume hypotheses \textrm{(H.1)-(H.4) }hold. Then, for any
fixed $k$ in $\mathbb{N},$ there exist $x_{k},y_{k}\in \Omega $ such that
the following estimates hold:%
\begin{equation}
\begin{array}{l}
\int_{\Omega }{\hat{u}_{R}}^{1+p(x)(d^{k}-1)}dx\leq \max \{1,|\Omega |\}\max
\{\Vert {\hat{u}_{R}}\Vert _{p(x)d^{k}}^{p(x_{k})d^{k}},\Vert {\hat{v}_{R}}%
\Vert _{q(x)d^{k}}^{q(y_{k})d^{k}}\},%
\end{array}
\label{36}
\end{equation}%
\begin{equation}
\begin{array}{l}
\int_{\Omega }{\hat{u}_{R}}|{\hat{u}_{R}}|^{\alpha (x)-1}|{\hat{v}_{R}}%
|^{\beta (x)+1}|{\hat{u}_{R}}|^{1+p(x)(d^{k}-1)}dx\leq 2\max \{\Vert {\hat{u}%
_{R}}\Vert _{p(x)d^{k}}^{p(x_{k})d^{k}},\Vert {\hat{v}_{R}}\Vert
_{q(x)d^{k}}^{q(y_{k})d^{k}}\},%
\end{array}
\label{37}
\end{equation}%
where $|\Omega |$ denotes the Lebesgue measure of a set $\Omega $ in $%
\mathbb{R}
^{N}$.
\end{lemma}

\textbf{Proof.} \ignorespaces
Before starting the proof, let us note that 
\begin{equation*}
\begin{array}{rl}
\frac{\alpha (x)+1+p(x)(d^{k}-1)}{p(x)d^{k}}+\frac{\beta (x)+1}{q(x)d^{k}}=
& \left[ \frac{\alpha (x)+1}{p(x)}+\frac{\beta (x)+1}{q(x)}\right] \frac{1}{%
d^{k}}+\frac{d^{k}-1}{d^{k}} \\ 
= & \frac{1}{d^{k}}+\frac{d^{k}-1}{d^{k}}=1,%
\end{array}
\label{holder}
\end{equation*}
%
%
%
where $d$ is chosen as in (\ref{53}).

Let us prove (\ref{36})\textbf{. }Since ${\hat{u}_{R}\in L}%
^{p(x)d^{k}}(\Omega )$ and $p(x)d^{k}>p(x)d^{k}-p(x)+1>0$ then ${\hat{u}%
_{R}\in L}^{\frac{p(x)d^{k}}{1+p(x)(d^{k}-1)}}(\Omega )$. Therefore, by H%
\"{o}lder's inequality and Mean value Theorem, there exist $x_{k}$ and $%
t_{k}\in \Omega $ such that%
\begin{equation*}
\begin{array}{l}
\int_{\Omega }{\hat{u}_{R}}^{1+p(x)(d^{k}-1)}dx\leq ||{\normalsize 1}%
_{\Omega }||_{d^{k}p^{\prime }(x)}\Vert {\hat{u}_{R}}\Vert
_{d^{k}p(x)}^{d^{k}p(x_{k})}=||{\normalsize 1}_{\Omega }||_{d^{k}p^{\prime
}(x)}^{\frac{d^{k}p^{\prime }(t_{k})}{d^{k}p^{\prime }(t_{k})}} \Vert {\hat{u%
}_{R}}\Vert _{d^{k}p(x)}^{d^{k}p(x_{k})} \\ 
=|\Omega |^{\frac{1}{d^{k}p^{\prime }(t_{k})}}\Vert {\hat{u}_{R}}\Vert
_{d^{k}p(x)}^{d^{k}p(x_{k})}\leq \max \{1,|\Omega |\}\Vert {\hat{u}_{R}}%
\Vert _{d^{k}p(x)}^{d^{k}p(x_{k})}.%
\end{array}%
\end{equation*}
This shows that the inequality (\ref{36}) holds true. Here $p^{\prime }$ and 
$p$ are conjugate variable exponents functions.

Next, we show\textbf{\ }(\ref{37})\textbf{.} By (\ref{holder}) and Young's
inequality, we get%
\begin{equation}
\begin{array}{l}
|\int_{\Omega }|{\hat{u}_{R}}|^{\alpha (x)+1+p(x)(d^{k}-1)}|{\hat{v}_{R}}%
|^{\beta (x)+1}dx|\leq \int_{\Omega }|{\hat{u}_{R}}|^{\alpha
(x)+1+p(x)(d^{k}-1)}|{\hat{v}_{R}}|^{\beta (x)+1}dx \\ 
\leq \int_{\Omega }\frac{\alpha (x)+1+p(x)(d^{k}-1)}{p(x)d^{k}}|{\hat{u}_{R}}%
|^{p(x)d^{k}} dx+\int_{\Omega }\frac{\beta (x)+1}{q(x)d^{k}}|{\hat{v}_{R}}%
|^{q(x)d^{k}} dx \\ 
\leq \int_{\Omega }|{\hat{u}_{R}}|^{p(x)d^{k}} dx+\int_{\Omega }|{\hat{v}_{R}%
}|^{q(x)d^{k}} dx.%
\end{array}
\label{45}
\end{equation}%
Observe from (\ref{12}) that%
\begin{equation*}
\begin{array}{l}
\int_{\Omega }|\frac{{\hat{u}_{R}}}{\Vert {\hat{u}_{R}}\Vert _{p(x)d^{k}}}%
|^{p(x)d^{k}} dx=1.%
\end{array}%
\end{equation*}
Using the mean value theorem, there exists $x_{k}\in \Omega $ such that%
\begin{equation}
\begin{array}{l}
\int_{\Omega }|{\hat{u}_{R}}|^{p(x)d^{k}} dx=\Vert {\hat{u}_{R}}\Vert
_{p(x)d^{k}}^{p(x_{k})d^{k}}.%
\end{array}
\label{46}
\end{equation}%
Similarly, we can find $y_{k}\in \Omega $ such that 
\begin{equation}
\begin{array}{l}
\int_{\Omega }|{\hat{v}_{R}}|^{q(x)d^{k}} dx=\Vert {\hat{v}_{R}}\Vert
_{q(x)d^{k}}^{q(y_{k})d^{k}}.%
\end{array}
\label{47}
\end{equation}%
Then, combining (\ref{45}), (\ref{46}) and (\ref{47}), the inequality (\ref%
{37}) holds true, ending the proof of the lemma \ref{estimk}. {\ \rule%
{0.2cm}{0.2cm}{\ }}

By using the Lemma we can prove the next result.

\begin{lemma}
\label{iteration} Assume \textrm{(H.1)-(H.4)} hold. Let $(\hat{u}_{R},\hat{v}%
_{R})\in X_{0}^{p(x),q(x)}(\Omega )$ be a solution of problem (\ref{p}).
Then, 
\begin{equation*}
({\hat{u}_{R}},{\hat{v}_{R}})\in L^{p(x)d^{k}}(\Omega )\times
L^{q(x)d^{k}}(\Omega ),\forall k\in \mathbb{N}.  \label{kk+1}
\end{equation*}
\end{lemma}

\textbf{Proof.} \ignorespaces
We employ a recursive reasoning. Since $(\hat{u}_{R},\hat{v}_{R})\in
X_{0}^{p(x),q(x)}(\Omega ),$ it is obvious that $({\hat{u}_{R}},{\hat{v}_{R}}%
)\in L^{p(x)}(\Omega )\times L^{q(x)}(\Omega ).$ So, (\ref{36}) remains true
for $k=0.$

{Assume that the conjecture \textquotedblleft $({\hat{u}_{R}},{\hat{v}_{R}}%
)\in L^{p(x)d^{l}}(\Omega )\times L^{q(x)d^{l}}(\Omega )$\textquotedblright
\thinspace \thinspace\ holds at every level $l\leq k$ and we claim that 
\begin{equation}
({\hat{u}_{R}},{\hat{v}_{R}})\in L^{p(x)d^{k+1}}(\Omega )\times
L^{q(x)d^{k+1}}(\Omega ).  \label{k+1}
\end{equation}%
To do it, we inserte $\varphi ={\hat{u}_{R}}^{1+p(x)(d^{k}-1)}$ in (\ref{1})
we get 
\begin{equation}
\begin{array}{l}
\int_{\Omega }\left\vert \nabla {\hat{u}_{R}}\right\vert ^{p(x)-2}\nabla {%
\hat{u}_{R}}\nabla {(\hat{u}_{R}}^{1+p(x)(d^{k}-1)}) dx \\ 
=\lambda _{R}^{\ast }\int_{\Omega }c(x)(\alpha (x)+1){\hat{u}_{R}}|{\hat{u}%
_{R}}|^{\alpha (x)-1}|{\hat{v}_{R}}|^{\beta (x)+1}{\hat{u}_{R}}%
^{1+p(x)(d^{k}-1)} dx.%
\end{array}
\label{39}
\end{equation}%
Observe that 
\begin{equation}
\begin{array}{l}
\int_{\Omega }\left\vert \nabla {\hat{u}_{R}}\right\vert ^{p(x)-2}\nabla {%
\hat{u}_{R}}\nabla \left( {\hat{u}_{R}}^{1+p(x)(d^{k}-1)}\right) dx \\ 
=\int_{\Omega }(d^{k}-1)\nabla {p}\nabla {\hat{u}_{R}}\left\vert \nabla {%
\hat{u}_{R}}\right\vert ^{p(x)-2}{\hat{u}_{R}}^{1+p(x)(d^{k}-1)}\ln {\hat{u}%
_{R}}\,\,dx \\ 
{\ \ \ }+\int_{\Omega }\left[ 1+p(x)(d^{k}-1)\right] \left\vert \nabla {\hat{%
u}_{R}}\right\vert ^{p(x)}{\hat{u}_{R}}^{p(x)(d^{k}-1)}dx.%
\end{array}
\label{40}
\end{equation}%
and 
\begin{equation}
\begin{array}{l}
\left\vert \nabla {\hat{u}_{R}}\right\vert ^{p(x)}{\hat{u}_{R}}%
^{p(x)(d^{k}-1)}=\frac{1}{d^{kp(x)}}|\nabla ({\hat{u}_{R}})^{d^{k}}|^{p(x)}.%
\end{array}%
\hat{u}_{R}  \label{41}
\end{equation}%
Then on the one hand%
\begin{equation}
\begin{array}{c}
\int_{\Omega }\frac{1+p(x)(d^{k}-1)}{d^{kp(x)}}|\nabla \left( {\hat{u}_{R}}%
\right) ^{d^{k}}|^{p\hat{u}_{R}(x)}dx\geq \int_{\Omega }\frac{d^{k}}{%
d^{kp(x)}}|\nabla \left( {\hat{u}_{R}}\right) ^{d^{k}}|^{p(x)}dx \\ 
\geq \frac{1}{d^{k(p^{+}-1)}}\int_{\Omega }|\nabla \left( {\hat{u}_{R}}%
\right) ^{d^{k}}|^{p(x)}dx,%
\end{array}
\label{42}
\end{equation}%
on the other hand, since ${\hat{u}_{R}}$ is assumed of class $C^{1}({%
\overline{\Omega }})$ and taking $\sup_{x\in \Omega }|\nabla
p|=M_{p}<+\infty ,$ we have 
\begin{equation}
\begin{array}{l}
\int_{\Omega }(d^{k}-1)\left\vert \nabla {p}\right\vert \left\vert \nabla {%
\hat{u}_{R}}\right\vert ^{p(x)-1}{\hat{u}_{R}}^{1+p(x)(d^{k}-1)}\left\vert
\ln {\hat{u}}\right\vert \,dx\leq {\hat{C}}M_{p}\int_{\Omega }{\hat{u}_{R}}%
^{1+p(x)(d^{k}-1)}\,dx,%
\end{array}
\label{43}
\end{equation}%
with some constant ${\hat{C}>0}$. Hence, gathering (\ref{39}), (\ref{40}), (%
\ref{42}) and (\ref{43}) together, one has%
\begin{equation}
\begin{array}{l}
\int_{\Omega }|\nabla \left( {\hat{u}_{R}}\right) ^{d^{k}}|^{p(x)}dx\leq
d^{k(p^{+}-1)}\int_{\Omega }\left[ 1+p(x)(d^{k}-1)\right] \left\vert \nabla {%
\hat{u}_{R}}\right\vert ^{p(x)}{\hat{u}_{R}}^{p(x)(d^{k}-1)}dx \\ 
\leq d^{k(p^{+}-1)}\int_{\Omega }(d^{k}-1)\left\vert \nabla {\hat{u}_{R}}%
\right\vert ^{p(x)-1}\left\vert \nabla {p}\right\vert {\hat{u}_{R}}%
^{1+p(x)(d^{k}-1)}\left\vert \ln {\hat{u}_{R}}\right\vert \,\,dx \\ 
\hbox{ \ \ \ \ \ \ \ \ \ \ \ \ \ }\,+\lambda _{R}^{\ast }\Vert c\Vert
_{\infty }(\alpha ^{+}+1)d^{k(p^{+}-1)}\int_{\Omega }|{\hat{u}_{R}}|^{\alpha
(x)+1+p(x)(d^{k}-1)}|{\hat{v}}|^{\beta (x)+1}dx \\ 
\leq {\hat{C}}_{p}d^{k(p^{+}-1)}\left[ \int_{\Omega }{\hat{u}_{R}}%
^{1+p(x)(d^{k}-1)}dx\,\,+\right. \\ 
\left. \hbox{ \ \ \ \ \ \ \ \ }\lambda _{R}^{\ast }\Vert c\Vert _{\infty
}(\alpha ^{+}+1)\int_{\Omega }|{\hat{u}_{R}}|^{\alpha (x)+1+p(x)(d^{k}-1)}|{%
\hat{v}_{R}}|^{\beta (x)+1}dx\right] ,%
\end{array}
\label{51}
\end{equation}%
where ${\hat{C}}_{p}=\max \{1,{\hat{C}}M_{p}\}$. }

Thanks to the use of the hypothesis \textrm{(H.3)}, the embeddings $L^{\pi
_{p}(x)}(\Omega )\hookrightarrow L^{dp(x)}(\Omega ),$ $W_{0}^{1,p(x)}(\Omega
)\hookrightarrow L^{\pi _{p}(x)}(\Omega )$ and $W_{0}^{1,p(x)}(\Omega
)\hookrightarrow L^{dp(x)}(\Omega )$ are continuous and thus, for any ${z}%
\in W_{0}^{1,p(x)}(\Omega ).$ We can conclude that there exists a constant $%
K>0$ so that%
\begin{equation}
\begin{array}{l}
\Vert {z}\Vert _{p(x)d}\leq K\Vert {z}\Vert _{1,p(x)}.%
\end{array}
\label{22}
\end{equation}%
From (\ref{12}) and through the mean value theorem observe that there exists 
$\xi _{k}\in \Omega $ such that 
\begin{equation*}
\begin{array}{l}
1=\int_{\Omega }\left\vert \frac{|\hat{u}_{R}|}{\Vert \hat{u}_{R}\Vert
_{p(x)d^{k+1}}}\right\vert ^{p(x)d^{k+1}}dx \\ 
=\int_{\Omega }\left\vert \frac{|\hat{u}_{R}|^{d^{k}}}{\Vert |\hat{u}%
_{R}|^{d^{k}}\Vert _{p(x)d}}\right\vert ^{p(x)d}\times \left( \frac{\Vert |%
\hat{u}_{R}|^{d^{k}}\Vert _{p(x)d}}{\Vert |\hat{u}_{R}|\Vert
_{p(x)d^{k+1}}^{d^{k}}}\right) ^{p(x)d}dx=\left( \frac{\Vert |\hat{u}%
_{R}|^{d^{k}}\Vert _{p(x)d}}{\Vert |\hat{u}_{R}|\Vert _{p(x)d^{k+1}}^{d^{k}}}%
\right) ^{p(\xi _{k})d},%
\end{array}%
\end{equation*}
which leads to%
\begin{equation}
\Vert |\hat{u}_{R}|^{d^{k}}\Vert _{p(x)d}=\Vert \hat{u}_{R}\Vert
_{p(x)d^{k+1}}^{d^{k}}.  \label{EQ}
\end{equation}%
Recalling from (\ref{normro}) that for every $z\in W_{0}^{1,p(x)}(\Omega
)\setminus \{0\}$%
\begin{equation}
\begin{array}{l}
\int_{\Omega }|\frac{|\nabla z|}{\Vert z\Vert _{1,p(x)}}|^{p(x)} dx=1.%
\end{array}
\label{23}
\end{equation}%
Applying (\ref{22}) and (\ref{23}) to $z={\hat{u}_{R}}^{d^{k}}$, besides the
mean value theorem and (\ref{EQ}), there exists $x_{k}\in \Omega $ such that%
\begin{equation}
\begin{array}{l}
K^{p(x_{k})}\int_{\Omega }|\nabla ({\hat{u}_{R})}%
^{d^{k}}|^{p(x)}dx=K^{p(x_{k})}\Vert {\hat{u}_{R}}^{d^{k}}\Vert
_{1,p(x)}^{p(x_{k})}=K^{p(x_{k})}\Vert {\hat{u}_{R}}^{d^{k}}\Vert
_{1,p(x)}^{p(x_{k})} \\ 
\geq \Vert {\hat{u}_{R}}^{d^{k}}\Vert _{dp(x)}^{p(x_{k})}=\Vert \hat{u}%
_{R}\Vert _{p(x)d^{k+1}}^{p(x_{k})d^{k}}=(\Vert \hat{u}_{R}\Vert
_{p(x)d^{k+1}}^{p(x_{k})d^{k+1}})^{\frac{1}{d}}.%
\end{array}
\label{52}
\end{equation}%
Combining (\ref{51}), (\ref{52}) with Lemma \ref{estimk}, we get the
following estimate 
\begin{equation}
\begin{array}{l}
\Vert {\hat{u}_{R}}\Vert _{p(x)d^{k+1}}^{p(x_{k})d^{k+1}}\leq
C_{1}d^{kd(p^{+}-1)}\left( \max \{\Vert {\hat{u}_{R}}\Vert
_{p(x)d^{k}}^{p(x_{k})d^{k}},\Vert {\hat{v}_{R}}\Vert
_{q(x)d^{k}}^{q(y_{k})d^{k}}\}\right) ^{d},%
\end{array}%
\end{equation}%
Acting also in (\ref{1}) with $\psi ={\hat{v}}^{1+q(x)(d^{k}-1)}$ and
repeating the argument above, we obtain%
\begin{equation}
\begin{array}{l}
\Vert {\hat{v}_{R}}\Vert _{q(x)d^{k+1}}^{q(x_{k})d^{k+1}}\leq
C_{2}d^{kd(q^{+}-1)}\left( \max \{\Vert {\hat{v}_{R}}\Vert
_{p(x)d^{k}}^{p(x_{k})d^{k}},\Vert {\hat{v}_{R}}\Vert
_{q(x)d^{k}}^{q(y_{k})d^{k}}\}\right) ^{d},%
\end{array}%
\end{equation}%
where $C_{1}$ and $C_{2}$ are two strictly positive constants. \newline
So, it derives 
\begin{equation}
\begin{array}{rl}
\max \{\Vert {\hat{u}_{R}}\Vert _{p(x)d^{k+1}}^{p(x_{k+1})d^{k+1}},\Vert {%
\hat{v}_{R}}\Vert _{q(x)d^{k+1}}^{q(y_{k+1})d^{k+1}}\}\leq & C_{3}d^{k\hat{d}%
}\left( \max \{\Vert {\hat{v}_{R}}\Vert _{p(x)d^{k}}^{p(x_{k})d^{k}},\Vert {%
\hat{v}_{R}}\Vert _{q(x)d^{k}}^{q(y_{k})d^{k}}\}\right) ^{d} \\ 
\leq & C_{3}d^{k\hat{d}}\left( \max \{\Vert {\hat{v}_{R}}\Vert
_{p(x)d^{k}}^{p(x_{k})d^{k}},\Vert {\hat{v}_{R}}\Vert
_{q(x)d^{k}}^{q(y_{k})d^{k}}\}\right) ^{d},%
\end{array}%
\end{equation}%
where $\hat{d}$ satisfies (\ref{50}) and $C_{3}=\max \{C_{1},C_{2}\}.$ 
\newline
\medskip \medskip

Before continuing, we distinguish the cases where $\Vert {\hat{u}_{R}}\Vert
_{p(x)d^{k+1}},$ $\Vert {\hat{v}_{R}}\Vert _{q(x)d^{k+1}},$ $\Vert {\hat{u}%
_{R}}\Vert _{p(x)d^{k}}$ and $\Vert {\hat{v}_{R}}\Vert _{p(x)d^{k}}$ are
each either less than one or either greater than one. Using \textrm{(H.4)}
and (\ref{53}) we obtain 
\begin{equation}
\begin{array}{rl}
\ln \left( \max \{\Vert {\hat{u}_{R}}\Vert _{p(x)d^{k+1}}^{d^{k+1}},\Vert {%
\hat{v}_{R}}\Vert _{q(x)d^{k+1}}^{d^{k+1}}\}\right) \leq & \!\!\!\ln
(C_{3}d^{k\hat{d}}) \\ 
& +\hat{d}\ln \left( \max \{\Vert {\hat{u}_{R}}\Vert
_{p(x)d^{k}}^{d^{k}},\Vert {\hat{v}_{R}}\Vert _{q(x)d^{k}}^{d^{k}}\}\right) .%
\end{array}
\label{lnuvk}
\end{equation}%
Now set%
\begin{equation}
\begin{array}{l}
E_{k}=\max \left\{ \ln \Vert {\hat{u}_{R}}\Vert _{p(x)d^{k}}^{d^{k}},\ln
\Vert {\hat{v}_{R}}\Vert _{q(x)d^{k}}^{d^{k}}\right\} \hbox{ \ and \ }\rho
_{k}=ak+b,%
\end{array}
\label{24}
\end{equation}%
with 
\begin{equation}
a=\ln d^{\hat{d}}\hbox{,\ \ \ \ }b=\ln C_{3}.
\end{equation}%
Then the recursive rule (\ref{lnuvk}) becomes 
\begin{equation}
\begin{array}{l}
E_{k+1}\leq \rho _{k}+\hat{d}E_{k},%
\end{array}
\label{28}
\end{equation}%
which in turn gives%
\begin{equation}
\begin{array}{l}
\begin{array}{l}
E_{k+1}\leq E\hat{d}^{k},%
\end{array}%
\end{array}
\label{32}
\end{equation}%
where%
\begin{equation}
\begin{array}{l}
E=E_{1}+\frac{b}{\hat{d}-1}+\frac{a\hat{d}}{(\hat{d}-1)^{2}}.%
\end{array}
\label{49}
\end{equation}%
Indeed, using (\ref{28}), (\ref{24}) and Lemma \ref{L6}, we get%
\begin{equation}
\begin{array}{l}
E_{k+1}\leq \rho _{k}+\hat{d}E_{k}\leq E_{k+1}\leq \rho _{k}+\hat{d}\rho
_{k-1}+\hat{d}^{2}E_{k-1} \\ 
\leq \rho _{k}+\hat{d}\rho _{k-1}+\hat{d}^{2}\rho _{k-2}+\hat{d}^{3}E_{k-2}
\\ 
... \\ 
\leq \sum_{i=0}^{k-1}\hat{d}^{i}\rho _{k-i}+\hat{d}^{k}E_{1}=\hat{d}%
^{k}(a\sum_{i=1}^{k}\frac{i}{\hat{d}^{i}}+b\sum_{i=1}^{k}\frac{1}{\hat{d}^{i}%
}+E_{1}) \\ 
\leq \hat{d}^{k}(\frac{a\hat{d}}{(\hat{d}-1)^{2}}+\frac{b}{\hat{d}-1}+E_{1})=%
\hat{d}^{k}E.%
\end{array}%
\end{equation}%
Here Lemma \ref{L6} is applied choosing $s=1/\hat{d}<1$ and $r=k+1$. So on,
according to (\ref{24}) and (\ref{32}), its follows that 
\begin{equation}
\begin{array}{l}
\max \{\Vert {\hat{u}_{R}}\Vert _{p(x)d^{k}},\Vert {\hat{v}_{R}}\Vert
_{q(x)d^{k}}\}\leq e^{E\frac{\max \{p^{+},q^{+}\}^{k-1}}{d}}.%
\end{array}%
\end{equation}%
We fix $k$ in $\mathbb{N},$ then we conclude that the assert (\ref{kk+1}) in
Lemma \ref{iteration} holds. The proof of Lemma \ref{iteration} is complete. 
{\ \rule{0.2cm}{0.2cm}{\ }}

\medskip\medskip

Now, let us end the proof of Theorem by showing that $({\hat{u}_R,\hat{v}_R)}
$ is bounded in $\Omega $.

\begin{lemma}
\label{linfini} Let $({\hat{u}_{R},\hat{v}_{R})}$ be a solution of (\ref{p})
corresponding to the eigenvalue $\lambda ^{\ast }.$ Assume that hypotheses 
\textrm{(H.1)-(H.4)} hold. Then, $\hat{u}_{R}$ and $\hat{v}_{R}$ are bounded
in $\Omega .$
\end{lemma}

\textbf{Proof.} \ignorespaces
Argue by contradiction. It means that we suppose that for all $L>0$, there
exists $\Omega _{L}\subset \Omega $, $|\Omega _{L}|>0$ such that for all $%
x\in \Omega _{L}$ we have $|{\hat{u}_{R}(x)|>L}$. Fix $k$ and choose $L$
large enough so that%
\begin{equation}
\begin{array}{l}
\frac{p^{-}\ln L}{p^{+}E\max \left\{ p^{+},q^{+}\right\} ^{k+1}}>1.%
\end{array}
\label{60}
\end{equation}%
From lemma \ref{L1} we get%
\begin{equation*}
\begin{array}{l}
L^{p^{-}d^{k+1}}|\Omega _{L}|\leq \int_{\Omega _{L}}L^{p(x)d^{k+1}}\,dx\leq
\int_{\Omega _{L}}|{\hat{u}_{R}|}^{p(x)d^{k+1}}\,dx \\ 
\leq \int_{\Omega }|{\hat{u}_{R}|}^{p(x)d^{k+1}}\,dx\leq \max \{\Vert {\hat{u%
}_{R}}\Vert _{p(x)d^{k+1}}^{p^{+}d^{k+1}},\Vert {\hat{u}_{R}}\Vert
_{p(x)d^{k+1}}^{p^{-}d^{k+1}}\}.%
\end{array}%
\end{equation*}
By (\ref{24}), (\ref{32}), and (\ref{53}) it follows that%
\begin{equation*}
\begin{array}{l}
d^{k+1}p^{-}\ln L+\ln |\Omega _{L}|\leq p^{+}E_{k+1}\leq p^{+}E\hat{d}^{k}%
\end{array}%
\end{equation*}
After using (\ref{60}) and dividing by $\hat{d}^{k+1},$ we get 
\begin{equation}
\begin{array}{l}
1+\frac{\ln |\Omega _{L}|}{\hat{d}^{k+1}}<1/\hat{d}.%
\end{array}
\label{L}
\end{equation}%
We choose $k$ sufficiently large in (\ref{L}). This forces $\hat{d}<1,$
which contradicts (\ref{50}). This proves the lemma \ref{linfini}. {\ \rule%
{0.2cm}{0.2cm}{\ }}

\medskip \medskip Next, we show that ${\hat{u}_{R}}$ and ${\hat{v}_{R}}$ are
strictly positive in $\Omega $.

\begin{lemma}
\label{positif} Let $({\hat{u}_{R},\hat{v}_{R})}$ be a solution of (\ref{p})
corresponding to the eigenvalue $\lambda ^{\ast }.$ Then, the following
asserts hold

\begin{enumerate}
\item $\hat{u}_{R}>0$ (resp. $\hat{v}_{R}>0$) in $\Omega .$

\item There exists $\delta \in (0,1)$ such that $\hat{u}_{R}$ is of class $%
C^{1,\delta }({\overline{\Omega }}).$
\end{enumerate}
\end{lemma}

\textbf{Proof.} \ignorespaces
\textbf{Step 1}.\textbf{\ }$\hat{u}_{R}\geq 0$\textbf{\ (resp. }$\hat{v}%
_{R}\geq 0$\textbf{\ in }$\Omega $\textbf{\ )} \newline
First, observe that%
\begin{equation*}
|u|=\max (u,0)+\min (u,0)\in W_{0}^{1,p(x)}(\Omega )
\end{equation*}%
and 
\begin{equation*}
|\nabla |u||\leq |\nabla \max (u,0)|+|\nabla \min (u,0)|\leq |\nabla u|.
\end{equation*}
Then it turns out that%
\begin{equation*}
\mathcal{A}(|{\hat{u}_{R}|},|{\hat{v}_{R}|})\leq \mathcal{A}({\hat{u}_{R}},{%
\hat{v}_{R}})\hbox{ and }\mathcal{B}(|{\hat{u}_{R}|},|{\hat{v}_{R}|})=%
\mathcal{B}({\hat{u}_{R}},{\hat{v}})=R.
\end{equation*}
Thereby (\ref{c1}) and (\ref{14}), it follows that%
\begin{equation*}
\mathcal{A}(|{\hat{u}_{R}|},|{\hat{v}_{R}|})\leq \mathcal{A}({\hat{u}_{R}},{%
\hat{v}_{R}})=R\lambda _{R}^{\ast }\leq \mathcal{A}(|{\hat{u}_{R}|},|{\hat{v}%
_{R}|}),
\end{equation*}
which implies that $\mathcal{A}(|{\hat{u}_{R}|},|{\hat{v}_{R}|})=R\lambda
_{R}^{\ast }$, showing that $(|{\hat{u}_{R}|},|{\hat{v}_{R}|})$ is a
solution of (\ref{p}). Therefore, we can assume that ${\hat{u}_{R}},{\hat{v}%
_{R}\geq 0}$ in $\Omega $. \newline
\medskip \medskip \newline
\textbf{Step 2. $\hat{u}_{R}>0$ (resp. $\hat{v}_{R}>0$) in $\Omega $ }%
\newline
Inspired by the ideas in \cite{LZZ}, let $m>0$ be a constant such that $%
h(\cdot )\in C^{2}(\overline{\partial \Omega }_{3m}),$ with $\overline{%
\partial \Omega }_{3m}=\{x\in \overline{\Omega }:h(x)\leq 3m\}$. Define the
functions%
\begin{equation*}
\mathcal{U}(x)=\left\{ 
\begin{array}{ll}
e^{\kappa h(x)}-1 & \hbox{ if }h(x)<\sigma _{1} \\ 
e^{\kappa h(x)}-1+\kappa e^{\kappa \sigma _{1}}\int_{\sigma _{1}}^{h(x)}(%
\frac{2m-t}{2m-\sigma _{1}})^{\frac{2}{p^{-}-1}}dt & \hbox{ if }\sigma
_{1}\leq h(x)<2\sigma _{1} \\ 
e^{\kappa h(x)}-1+\kappa e^{\sigma _{1}}\int_{\sigma _{1}}^{2m}(\frac{2m-t}{%
2m-\sigma _{1}})^{\frac{2}{p^{-}-1}}dt & \hbox{ if }2\sigma _{1}\leq h(x)%
\end{array}%
\right.
\end{equation*}
and%
\begin{equation*}
\mathcal{V}(x)=\left\{ 
\begin{array}{ll}
e^{\kappa h(x)}-1 & \hbox{ if }h(x)<\sigma _{2} \\ 
e^{\kappa h(x)}-1+\kappa e^{\kappa \sigma _{2}}\int_{\sigma _{2}}^{h(x)}(%
\frac{2m-t}{2m-\sigma _{2}})^{\frac{2}{q^{-}-1}}dt & \hbox{ if }\sigma
_{2}\leq h(x)<2\sigma _{2} \\ 
e^{\kappa h(x)}-1+\kappa e^{\kappa \sigma _{2}}\int_{\sigma _{2}}^{2m}(\frac{%
2m-t}{2m-\sigma _{2}})^{\frac{2}{q^{-}-1}}dt & \hbox{ if }2\sigma _{2}\leq
h(x),%
\end{array}%
\right.
\end{equation*}
where $(\sigma _{1},\sigma _{2})=(\frac{\ln 2}{\kappa p^{+}},\frac{\ln 2}{%
\kappa q^{+}})$ and $\kappa >0$ is a parameter. A quite similar calculations
as in \cite[pages 11 and 12]{LZZ} furnish 
\begin{equation}
\begin{array}{l}
-\Delta _{p(x)}({\mu }_{1}\mathcal{U)}\leq \lambda _{R}^{\ast }c(x)(\alpha
(x)+1)(\mu _{1}\mathcal{U})^{\alpha (x)}{\hat{v}_{R}}^{\beta (x)+1}%
\hbox{ \
in }\Omega%
\end{array}
\label{lmu1}
\end{equation}%
and%
\begin{equation}
\begin{array}{l}
-\Delta _{q(x)}{(\mu }_{2}{\mathcal{V})}\leq \lambda _{R}^{\ast }c(x)(\beta
(x)+1){\hat{u}_{R}}^{\alpha (x)+1}(\mu _{2}\mathcal{V}{)}^{\beta (x)}%
\hbox{
\ in }\Omega ,%
\end{array}
\label{lmu2}
\end{equation}%
where $\mu _{1}=\exp (\kappa \frac{1-p^{-}}{\max_{\overline{\Omega }}|\nabla
p|+1})\ \ $and $\ \mu _{2}=\exp (\kappa \frac{1-q^{-}}{\max_{\overline{%
\Omega }}|\nabla q|+1})$, provided that $\kappa >0$ is large enough. \newline
Now, for any $(z,w)\in X_{0}^{p(x),q(x)}(\Omega ),$ denote by 
\begin{equation*}
\mathcal{L}_{p}(z,w)=-\Delta _{p(x)}z-\lambda _{R}^{\ast }c(x)(\alpha
(x)+1)z|z|^{\alpha (x)-1}|w|^{\beta (x)+1}
\end{equation*}
and 
\begin{equation*}
\mathcal{L}_{q}(z,w)=-\Delta _{q(x)}w-\lambda _{R}^{\ast }c(x)(\beta
(x)+1)|z|^{\alpha (x)+1}w|w|^{\beta (x)-1},
\end{equation*}
(\ref{lmu1}) and (\ref{lmu2}) may be formulated respectively as follows 
\begin{equation*}
\mathcal{L}_{p}(\mu _{1}\mathcal{U},\hat{v}_{R})\leq 0\,\,\hbox{ and }%
\mathcal{L}_{q}(\hat{u}_{R},\mu _{2}\mathcal{V})\leq 0,\,\,\hbox{ in }\Omega.
\end{equation*}
Hence, from the above notation, we get%
\begin{equation*}
\begin{array}{l}
\mathcal{L}_{p}(\mu _{1}\mathcal{U},\hat{v}_{R})\leq 0\leq \mathcal{L}_{p}(%
\hat{u}_{R},\hat{v}_{R})\hbox{ \ in }\Omega%
\end{array}%
\end{equation*}
and%
\begin{equation*}
\begin{array}{l}
\mathcal{L}_{q}(\hat{u}_{R},\mu _{2}\mathcal{V})\leq 0\leq \mathcal{L}_{q}(%
\hat{u}_{R},\hat{v}_{R})\hbox{ \ in }\Omega .%
\end{array}%
\end{equation*}
Since ${\mu }_{1}\mathcal{U}=\hat{u}_{R}=0$ and ${\mu }_{2}{\mathcal{V}}=%
\hat{v}_{R}=0$ on $\partial \Omega $, we are allowed to apply \cite[Lemma 2.3%
]{QZ} and we deduce that%
\begin{equation*}
\begin{array}{l}
{\hat{u}_{R}\geq \mu }_{1}\mathcal{U}>0\hbox{ \ and \ }{\hat{v}_{R}\geq {\mu 
}_{2}\mathcal{V}>0}\hbox{ in }\Omega.%
\end{array}%
\end{equation*}
Thereby the positivity of $({\hat{u}_{R}},{\hat{v}_{R}})$ in $\Omega $ is
proven. \newline
To end the proof of Lemma \ref{linfini}, we claim a regularity property for $%
{\hat{u}_{R}}$ and ${\hat{v}_{R}}.$ \medskip \medskip

\textbf{Step 3. Regularity property}\medskip \newline
For $p,q\in C^{1}(\overline{\Omega })\cap C^{0,\theta }(\overline{\Omega })$
for certain $\theta \in (0,1)$, owing to \cite[Theorem 1.2]{Fan} the
solution $({\hat{u}_{R}},{\hat{v}_{R}})$ belongs to $C^{1,\delta }(\overline{%
\Omega })\times C^{1,\delta }(\overline{\Omega })$ for certain $\delta \in
(0,1)$. This completes the proof. {\ \rule{0.2cm}{0.2cm}{\ }}

\end{document}